\documentclass[12pt]{article}
\usepackage{amsfonts}
\usepackage{amssymb}
\usepackage{latexsym}

\usepackage{epsf,graphicx
}

\newcommand{\pl}{\partial}
\newcommand{\pr}{\noindent{\bf Proof.}\quad }

\newcommand{\be} {\begin{eqnarray}}
\newcommand{\ee} {\end{eqnarray}}
\newcommand{\bep} {\begin{eqnarray*}}
\newcommand{\eep} {\end{eqnarray*}}
\newcommand{\dst} {\displaystyle}

\newcommand {\spi}{\mathop{\rm Sp}\nolimits}
\newcommand {\ud}{\mathop{\rm Univ}\nolimits(\Delta)}
\newcommand {\Hol}{\mathop{\rm Hol}\nolimits}
\renewcommand {\Im}{\mathop{\rm Im}\nolimits}
\renewcommand {\Re}{\mathop{\rm Re}\nolimits}
\newcommand {\bt}{\mathop{\mathcal{B}[\tau]}\nolimits}
\newcommand {\bu}{\mathop{\mathcal{B}[1]}\nolimits}
\newcommand {\gt}{\mathop{\mathcal{G}[\tau]}\nolimits}
\newcommand {\gu}{\mathop{\mathcal{G}[1]}\nolimits}
\newcommand {\gtp}{\mathop{\mathcal{G}^+[\tau]}\nolimits}
\newcommand {\gup}{\mathop{\mathcal{G}^+[1]}\nolimits}
\newcommand {\Pp}{\mathop{\mathcal P}\nolimits}
\newcommand {\pup}{\mathop{\mathcal P}^+[1]}
\newcommand {\Ss}{\mathcal{S}}
\newcommand {\G}{\mathop{\mathcal G}\nolimits}

\newfont{\bbb}{msbm10 at 12pt}
\newfont{\bbbb}{msbm10 at 8pt}
\def\Bbb#1{\hbox{\bbb #1}}
\def\Bbbb#1{\hbox{\bbbb #1}}
\newcommand{\R}{{\Bbb R}}

\newcommand{\N}{{\Bbb N}}
\newcommand{\NN}{{\Bbbb N}}

\newcommand{\C}{{\Bbb C}}

\newtheorem{defin}{Definition}[section]
\newtheorem{theorem}{Theorem}[section]
\newtheorem{corol}{Corollary}[section]

\newcounter{rem}
\newcommand{\rema}{{\bf Remark \addtocounter{rem}{1}\therem.}}

\begin{document}

\title{Controlled approximation and interpolation for some classes of holomorphic functions}

\author{Mark Elin \\ {\small Department of Mathematics, ORT  Braude College,} \\
{\small P.O. Box 78, Karmiel 21982, ISRAEL} \\ {\small e-mail:
mark.elin@gmail.com}
\\ David Shoikhet \\ {\small Department of Mathematics, ORT Braude College,}
\\ {\small P.O. Box 78, Karmiel 21982, ISRAEL} \\ {\small e-mail: davs27@netvision.net.il}
\\ Lawrence Zalcman\footnote{Research supported by the German--Israel Foundation for Scientific
Research and Development, G.I.F. Grant No. I-809-234-6/2003.}
\\ {\small Department of Mathematics, Bar-Ilan University,} \\
{\small 52900 Ramat-Gan, ISRAEL} \\ {\small e-mail:
zalcman@macs.biu.ac.il}  }
\date{}


\maketitle

\begin{abstract}
This paper reports on constructive approximation methods for three
classes of holomorphic functions on the unit disk which are
closely connected each other: the class of starlike and spirallike
functions, the class of semigroup generators, and the class of
functions with positive real part. It is more-or-less known that
starlike or spirallike functions can be defined as solutions of
singular differential equations which, in general, are not stable
under the motion of interior singular points to the boundary. At
the same time, one can establish a perturbation formula which
continuously transforms a starlike (or spirallike) function with
respect to a boundary point to a starlike (or spirallike) function
with respect to an interior point. This formula is based on an
appropriate approximation method of holomorphic generators which
determine the above-mentioned differential equations. In turn, the
well-known Berkson--Porta parametric representation of holomorphic
generators leads us to study an approximation-interpolation
problem for the class of functions with positive real part. While
this problem is of independent interest, the solution we present
here is again based on the Berkson--Porta formula. Finally, we
apply our results to solve a natural perturbation problem for
one-parameter semigroups of holomorphic self-mappings, as well as
the eigenvalue problem for the semigroup of composition operators.
\end{abstract}

\section{Preliminaries}

\setcounter{equation}{0}

\subsection{Starlike and spirallike functions}
Let $\Delta=\{z\in\C:\,|z|<1\}$ be the open unit disk in the
complex plane $\C$, and $\Pi=\{z\in\C:\,\Re z\ge0\}$ the right
half-plane. We denote the set of holomorphic functions on $\Delta$
which take values in a set $\Omega\subset\C$ by
$\Hol(\Delta,\Omega)$. In particular, $\Hol(\Delta,\C)$ is the set
of all holomorphic functions on $\Delta$. This set is a Frech\'et
space endowed with the seminorms $\|f\|_K:=\max_{z\in K}|f(z)|,$
where $K$ is a compact subset of $\Delta$.

For brevity, we write $\Hol(\Delta)$ for $\Hol(\Delta,\Delta)$ and
$\Pp$ for $\Hol(\Delta,\Pi)$. The set $\Hol(\Delta)$ is a
semigroup with respect to composition.

We denote the subset of $\Hol(\Delta,\C)$ consisting of univalent
functions (without any normalization) by $\ud$.

\begin{defin}
A univalent function $h$ is said to be {\textsf{starlike} } if for
each point $z\in\Delta$, the linear segment joining $h(z)$ to zero
\[
(0,h(z)]:=\{th(z):\ t\in(0,1]\}
\]
lies entirely in the image $h(\Delta)$.
\end{defin}

\begin{defin}
A univalent function $h$ is said to be {\sf spirallike} if there
is a number $\mu\in\C$ with $\Re\mu>0$ such that for every point
$z\in\Delta$, the spiral curve
\[
\{e^{-t\mu}h(z):\ t\ge 0\}
\]
lies entirely in $h(\Delta)$.
\end{defin}

It is clear that if $\mu$ can be chosen to be real, then the
spirallike function $h$ is actually starlike.

It follows from the above definitions that
$0\in\overline{h(\Delta)}$ for each starlike or spirallike
function $h$. We distinguish two situations.

\noindent$\bullet$ In case $0\in h(\Delta)$, the starlike
(spirallike) function $h$ is said to be starlike (spirallike) with
respect to an interior point.

\noindent$\bullet$ In case $0\in\pl h(\Delta)$, the function $h$
is said to be starlike (spirallike) with respect to a boundary
point.

For any starlike (spirallike) function with respect to an interior
point, there exists a unique point $\tau\in\Delta$ such that
$h(\tau)=0$.

If $h$ is starlike (spirallike) with respect to a boundary point,
then one can show (see, for example, \cite{E-R-S2000, E-R-S2001a})
that there exists a unique point $\tau\in\pl\Delta$ such that the
value $h(\tau)$ defined as the angular limit at this point is
equal to zero: $h(\tau):=\angle\lim\limits_{z\to\tau}h(z)=0$.

In both cases, we use the notation $S^*[\tau]$ for the set of
starlike functions (with respect to an interior or a boundary
point) satisfying the condition $h(\tau)=0,\,
\tau\in\overline\Delta$, and the notation $\spi[\tau]$ for the set
of spirallike functions satisfying the same condition.

In this paper, we use autonomic dynamical systems to study
approximation problems for starlike and spirallike functions with
respect to a boundary point and for some related classes of
functions. Since such systems are time-independent, their
solutions form one-parameter semigroups of holomorphic
self-mappings of the open unit disk.

\vspace{3mm}

\subsection{Semigroups of holomorphic self-mappings}
\begin{defin}
A family $\Ss=\{F_t\}_{t\ge0}\subset\Hol(\Delta)$ of holomorphic
self-mappings of $\Delta$ is called a {\sf one-parameter
continuous semigroup} (respectively, {\sf group}) if

(i) $F_{t+s}=F_t\circ F_s$ whenever $s,t$ and $s+t$ belong to
$\R^+$ (respectively, $\R$);

(ii) $F_0(z)=z$ for all $z\in\Delta$, that is, $F_0$ is  the
identity mapping on $\Delta$.

(iii) $\lim\limits_{t\to s}F_t(z)=F_s(z)$ for all $t>0,\,s\ge0$
(respectively, $t,s\in\R$) and for all $z\in\Delta$.
\end{defin}

The following result is due to Berkson and Porta \cite{BE-PH} (see
also \cite{RS-SD-96, RS-SD-98, AM-92}).

\begin{theorem}[{\rm see Propositions 3.2.1 and 3.2.2 \cite{SD}}]\label{2.1}
Let $\Ss=\left\{F_t\right\}_{t\ge0}$ be a one-parameter semigroup
of holomorphic self-mappings of $\Delta$ such that for each
$z\in\Delta$,
\begin{equation}
\lim_{t\rightarrow 0^{+}}F_t(z)=z.
\end{equation}
Then for each $z\in \Delta$, the limit
\begin{equation}\label{gen}
\lim_{t\rightarrow 0^{+}}\frac{z-F_t(z)}t=f(z),
\end{equation}
exists and is a holomorphic function on $\Delta$. The convergence
in (\ref{gen}) is uniform on each subset strictly inside $\Delta$.
Moreover, the semigroup $\Ss$ can be defined as the (unique)
solution of the Cauchy problem
\be\label{cauchy}
\left\{
\begin{array}{l}
{\displaystyle\frac{\partial F_t(z)}{\partial t}}+
f(F_t(z))=0,\quad t\ge0, \vspace{3mm}\\ F_0(z)=z,\quad z\in
\Delta.
\end{array}
\right.
\ee
\end{theorem}

\begin{defin}
Let $\Ss=\left\{F_t\right\}_{t\ge0}$ be a one-parameter continuous
semigroup of holomorphic self-mappings of $\Delta$. The function
$f\in\Hol(\Delta,\C)$ defined by the limit (\ref{gen}) is called
the {\sf (infinitesimal) generator} of $\Ss$.
\end{defin}

We denote the family of all holomorphic generators on $\Delta$ by
$\G$. This set is a real cone in $\Hol(\Delta,\C)$ \cite{RS-SD-96,
SD}.

A continuous version of the Denjoy--Wolff Theorem (see
\cite{RS-SD-97b}) asserts that {\it if a semigroup
$\Ss=\left\{F_t\right\}_{t\ge0}$ does not contain an elliptic
automorphism, then there is a unique point
$\tau\in\overline\Delta$ such that
$\lim\limits_{t\to\infty}F_t(z)=\tau$ for all $z\in\Delta$.}

This point $\tau$ is called the Denjoy--Wolff point of the
semigroup $\Ss$. For a given $\tau\in\overline\Delta$, the set of
all semigroups for which $\tau$ is their Denjoy--Wolff point is
denoted by $\bt$. The set of functions $f\in\G$ such that the
semigroup generated by $f$ belongs to $\bt$ is denoted by $\gt$.
Note that $\gt$ is a real subcone of $\G$.

\noindent$\bullet$ If $\Ss=\left\{F_t\right\}_{t\ge0}\in\bt$,
where $\tau\in\Delta$ is an interior point of $\Delta$, then
$\tau$ must be the unique (interior) fixed point of the semigroup
$\Ss$, i.e., $F_t(\tau)=\tau$ for all $t\ge0$.

Also, it follows by the uniqueness of the solution of the Cauchy
problem (\ref{cauchy}) that, in this case, $\tau$ is the unique
zero of $f$ in $\Delta$. Moreover, it can be shown (see, for
example, Theorem \ref{ThC} below) that $f\in\G$ belongs to $\gt,\
\tau\in\Delta,$ if and only if $f(\tau)=0$ and $\Re f'(\tau)>0$.

\noindent$\bullet$ If $\Ss=\left\{F_t\right\}_{t\ge0}\in\bt$,
where $\tau\in\pl\Delta$ is a boundary point, then its generator
$f$ does not vanish inside $\Delta$, but
\be\label{f0}
\angle\lim_{z\to\tau}f(z)=0
\ee
and
\[
\angle\lim_{z\to\tau}F_t(z)=\tau
\]
(see \cite{E-S2001}), where the symbol $\angle\lim$ denotes the
so-called angular (or non-tangential) limit at a boundary point of
$\Delta$ (see, for example, \cite{PC-92}).

Note, however, that $f\in\G$ may have more than one boundary null
point in the sense of (\ref{f0}). Moreover, in general, condition
(\ref{f0}) does not even imply that $\tau$ is a fixed point for
the semigroup $\Ss=\left\{F_t\right\}_{t\ge0}$ generated by $f$
(consider, for example, $f\in\G$ defined by $f(z)=z\sqrt{1-z}$
with $\tau=1$).

A characterization of the classes $\gt$, as well as the whole
class $\G$, is given in the following assertion.

\begin{theorem}[see \cite{BE-PH, A-E-R-S, SD}]\label{ThC}
The following statements are equivalent:

(i) $f\in\G$;

(ii) $f(z)=a-\bar{a}z^2+zp(z)$ for some $a\in\C$ and $p\in\Pp$;

(iii) $f$ admits the representation
\be\label{bp}
f(z)=(z-\tau)(1-z\bar\tau)p(z)
\ee
for some $\tau\in\overline\Delta$ and $p\in\Pp$.

Moreover, for a given $\tau\in\overline\Delta$, every generator
$f\in\gt$ admits the representation (\ref{bp}). Conversely, if $f$
is represented by (\ref{bp}) and does not generate a group of
elliptic automorphisms, then it belongs to $\gt$.
\end{theorem}
The parametric representation (\ref{bp}) is due to Berkson and
Porta \cite{BE-PH}. This representation is unique for
$\tau\in\overline\Delta$ and $p\in\Pp$.

The next theorem characterizes the class $\gt$ specifically for
the case where $\tau\in\pl\Delta$ is a boundary Denjoy--Wolff
point of generated semigroups.

\begin{theorem}[see \cite{E-S2001}]\label{ThD}
Let $f\in\G$ be a semigroup generator on $\Delta$ and let
$\Ss=\{F_t\}_{t\ge 0}$ be the semigroup generated by $f$. The
following are equivalent.

(i) $f$ has no null point in $\Delta$.

(ii) There is a point $\tau\in\pl\Delta$ such that
\[
\beta=\angle\lim\limits_{z\rightarrow\tau}{\frac{f(z)}{z-\tau}}=:\angle f'(\tau)
\]
exists finitely and $\Re{\beta}\ge 0$. Consequently, $f(\tau)=0$.

(iii) There exists a point $\tau\in\pl\Delta$ and a real
nonnegative number $\gamma$ such that
\[
\frac{|F_t(z)-\tau|^2}{1-|F_t(z)|^2} \le
\exp{(-t\gamma)}\,\frac{|z-\tau|^2}{1-|z|^2}
\]
for all $z\in\Delta$. Consequently, $f\in\gt$.

Moreover,

(a) the boundary points $\tau $ in (ii)--(iii) are the same;

(b) the limit $\beta$ in (ii) is a real nonnegative number;

(c) the maximal number $\gamma\ge 0$ for which (iii) holds
coincides with the limit value $\beta$ in~(ii).
\end{theorem}

This theorem can be considered the infinitesimal version of the
Julia--Wolff--Carath\'eodory Theorem. Assertion (iii) means that
all horocycles internally tangent to the boundary $\pl\Delta$ at
the point $\tau$ are invariant under the semigroup
action\footnote{In this case, the point $\tau\in\pl\Delta$ is
sometimes also called the sink point of the semigroup.}. Moreover,
if $\beta=f'(\tau)>0$, then choosing $\gamma=\beta$, we get an
exponential rate of convergence of the semigroup
$\Ss=\{F_t\}_{t\ge 0}$ to~$\tau$\vspace{2mm}.  Since we are mostly
interested in such a situation, we denote by $\gtp$
($\tau\in\pl\Delta$) the set of all generators satisfying the
conditions
\[
\angle\lim_{z\to\tau}f(z)=0,\quad
\angle\lim_{z\to\tau}\frac{f(z)}{z-\tau}=\beta>0.
\]

\vspace{3mm}

\subsection{Differential and functional equations for starlike and
spirallike functions} The dynamic approach to the study of
starlike and spirallike functions is based on the following
observation. Let $h\in\ud$; then $h$ is spirallike (starlike) if
and only if there exists $\mu\in\C$ with $\Re\mu>0$
($\mu\in(0,\infty)$) such that for each $t\ge0$, the function
\be\label{semig}
F_t(z)=h^{-1}\left(e^{-\mu t} h(z) \right)
\ee
is a well-defined holomorphic self-mapping of the open unit disk
$\Delta$.

Obviously, the family
$\Ss=\left\{F_t\right\}_{t\ge0}\subset\Hol(\Delta)$\vspace{2mm}
forms a one-parameter continuous semigroup. Note also that
$h\in\spi[\tau]$ (or $h\in S^*[\tau]$) if and only if $\Ss\in\bt$.

Now, differentiating (\ref{semig}) at $t=0^+$, we see that the
function $h$ satisfies the differential equation
\be\label{hgen}
\mu h(z)=h'(z)f(z),
\ee
where $f=\lim\limits_{t\to0^+}\frac{z-F_t(z)}t$ is the generator
of the semigroup $\Ss=\left\{F_t\right\}_{t\ge0}$.

In fact, the converse assertion also holds. However, as there are
some differences between the interior or the boundary location of
the the point $\tau$, we describe these cases separately.

\begin{theorem}\label{ThA}
Let $\tau\in\Delta$. A function $h\in\Hol(\Delta,\C)$ belongs to
$\spi[\tau]$ if and only if it is locally univalent and satisfies
equation (\ref{hgen}) for some $f\in\gt$ and $\mu=f'(\tau)$.
Moreover, in this case, $h$ also satisfies equation (\ref{semig}),
where $\Ss=\left\{F_t\right\}_{t\ge0}$ is the semigroup generated
by $f$.
\end{theorem}

\begin{theorem}\label{ThB}
Let $\tau\in\pl\Delta$. A univalent function $h\in\ud$ belongs to
$\spi[\tau]$ if and only if it satisfies equation (\ref{hgen}) for
some $f\in\gtp$ and $\mu\in\C$ with $\Re\mu>0$.
\end{theorem}

\noindent{\rema}  Note that the set $\G$ of holomorphic generators
is a real cone \cite{RS-SD-97b, SD}. Hence, for each spirallike
(starlike) function $h$, the pair $(\mu,f)$ in equation
(\ref{hgen}) can be replaced by a pair $(\alpha\mu,\alpha f)$ for
any $\alpha>0$. In fact, a function $h\in\ud$ satisfying equation
(\ref{hgen}) with some $\mu\in\C,\ \Re\mu>0,$ and $f\in\G$ is
actually starlike if and only if $\frac1\mu\cdot f\in\G$.

\vspace{2mm}

\noindent{\rema} Note that for $\tau\in\Delta$, equation
(\ref{hgen}) has no holomorphic solution when $\mu\not=f'(\tau)$,
since this equation has an interior singular point. In other
words, the number $\mu\in\C$ is uniquely determined by $f\in\gt$.

For $\tau\in\pl\Delta$, equation (\ref{hgen}) has a holomorphic
solution for each $\mu\in\C$. However, it can be shown (see
Theorem \ref{t3.8} below) that this solution is univalent if and
only if $\mu$ lies in the set $\Omega=\Omega_+\bigcup\Omega_-$,
where\linebreak $\Omega_\pm=\biggl\{w\in\C:\
|w\mp\beta|\le\beta,\,w\not=0\biggr\},\ \beta=f'(\tau)>0$.
Moreover, the solution $h$ of (\ref{hgen}) normalized by the
condition $h(0)=1$ belongs to $\spi[\tau]$ (respectively, to
$S^*[\tau]$) if and only if $\mu\in\Omega_+$ (respectively,
$0<\mu\le2\beta$).

\vspace{2mm}

\noindent{\rema} Both Theorem~\ref{ThA} and \ref{ThB} state that
for a given semigroup $\Ss=\left\{F_t\right\}_{t\ge0}$ generated
by $f\in\gt,\ \tau\in\overline\Delta,$ the solution $h$ of the
differential equation (\ref{hgen}) is a solution of the so-called
Schr\"oder equation
\be\label{schroder1}
h(F_t(z))=\lambda_t h(z)
\ee
with $\lambda_t=e^{-\mu t}$ \cite{CCC-MBD}.

On the other hand, if we do not require the univalence of
solutions, equation (\ref{schroder1}) can be considered the
eigenvalue problem for the semigroup of composition operators
$\left\{C_t\right\}_{t\ge0}$ defined by
\[
C_t(g)\biggl(=C_{F_t}(g)\biggr):=g(F_t)\quad\mbox{for all
}g\in\Hol(\Delta,\C).
\]

For $\tau\in\Delta$, the spectrum $\sigma(C_{F_t})$ is discrete,
while for $\tau\in\pl\Delta$, the spectrum $\sigma(C_{F_t})$
covers the whole region (see the discussion in Section 3 below).

\vspace{2mm}

In general, to each $f\in\Hol(\Delta,\C)$ there corresponds a
vector field $\Gamma_f$ defined by
\begin{equation}\label{gf}
\Gamma_f(g)(z):=g'(z)f(z)\quad\mbox{for all }g\in\Hol(\Delta,\C).
\end{equation}

Each vector field $\Gamma_f$ is locally integrable in the sense
that for each $z\in\Delta$, there exists a neighborhood $U$ of $z$
and a number $\delta>0$ such that the Cauchy problem
\be\label{cauchy1}
\left\{
\begin{array}{l}
\dst\frac{\partial u(t,z)}{\partial t}+f(u(t,z))=0 \vspace{3mm}\\
u(0,z) =z
\end{array}
\right.
\ee
has a unique solution $\{u(t,z)\}\subset\Delta$ defined on the set
$\{|t|<\delta\}\times U\subset\R\times\Delta.$

\begin{defin}[see \cite{AJ, UH, RS-SD-96}]
A vector field $\Gamma_f$ is said to be {\sf semi-complete}
(respectively, {\sf complete}) on $\Delta$ if the solution of the
Cauchy problem (\ref{cauchy1}) is well-defined on all of
$\R^+\times\Delta$ (respectively, $\R\times\Delta$).
\end{defin}

\noindent{\rema} In fact, Theorem~\ref{2.1} asserts that {\it a
vector field $\Gamma_f$ on the (Frech\'et) space $\Hol(\Delta,\C)$
is semi-complete if and only if the function $f$ is a generator of
a one-parameter continuous semigroup of holomorphic self-mappings
of $\Delta$.} In this case, the linear semigroup of composition
operators $\left\{C_t\right\}_{t\ge0}=
\left\{C_{F_t}\right\}_{t\ge0}$ is differentiable; and the
semi-complete vector field $\Gamma_f$ defined by (\ref{gf}) is its
generator in the space $\Hol(\Delta,\C)$,~i.e.,
\[
\lim_{t\to0^+}\frac1t\left(I-C_t\right)(g)= \Gamma_f(g)
\]
for each $g\in\Hol(\Delta,\C)$, where $I$ denotes the identity
operator on $\Hol(\Delta,\C)$.

Furthermore, the differential equation (\ref{hgen}), rewritten as
\be\label{hgen1}
\mu h=\Gamma_fh,
\ee
can be considered the eigenvalue problem for the linear operator
$\Gamma_f$.

\vspace{3mm}

\section{Approximation problems}
\setcounter{equation}{0}

\subsection{A perturbation problem for spirallike function}
The first approximation problem we discuss here can be described
as follows. Let $f\in\gup$, and let $h\in\spi[1]$ be a spirallike
function which satisfies the equation (\ref{hgen}):
\[
\mu h(z)=h'(z)f(z).
\]
Consider the perturbed equation
\[
\mu_\tau  h_\tau (z)= h_\tau '(z)f_\tau (z),
\]
where $f_\tau\in\gt, \ \tau\in\Delta$, and
$\mu_\tau=f_\tau'(\tau),$ are such that $f_\tau$ converges to $f$
locally uniformly on $\Delta$ when $\tau$ goes to $1$
unrestrictedly. (Note that such a perturbation is always possible,
see, for example, Theorem \ref{ThC}.) We ask: {\bf does the net
$\left\{h_\tau\right\}$ converge to $h$ as $\tau\to1$?}

The following simple example shows that, in general, the answer is
negative.

\vspace{2mm}

\noindent{\bf Example 1.} Let $f(z)=(z-1)\in\gup$. Then the
function $h(z)={1-z}$ satisfies equation (\ref{hgen}) with
$\mu=1$:
\[
h(z)=h'(z)f(z).
\]
Define now $\dst f_\tau(z)=\frac{(z-\tau)(1-z\tau)}{1-z},\
\tau\in(0,1)$. Obviously, $f_\tau\in\gt$ (by Theorem \ref{ThC})
and $f_\tau$ converges to $f$ as $\tau\to1^-$. Consider the
perturbed problem
\be\label{aux1}
\mu_\tau h_\tau (z)= h'_\tau (z)f_\tau(z),
\ee
where $\mu_\tau=f'_\tau(\tau)=1+\tau,\ f_\tau\to f$ as
$\tau\to1^-$. Then the function
\[
h_\tau(z)= \frac{(\tau-z)(1-z\tau)^\frac1\tau}{\tau}
\]
is a solution of equation (\ref{aux1}) satisfying
$h_\tau(0)=h(0)=1$. Letting $\tau$ tend to the boundary point $1$,
we obtain that the limit function
\[
\lim _{\tau\rightarrow 1 ^-}h_\tau(z)=(1-z)^2,
\]
which is different from $h(z)$.

At the same time if we choose $f_\tau\in\gt$ in a different way,
say $f_\tau(z)=z-\tau$, we see that $h_\tau=1-\frac z\tau$ defined
as a solution of equation (\ref{aux1}) converges to $h$.

\vspace{2mm}

Thus one can consider the following perturbation problem. {\bf For
any $\tau\in\Delta$, find a perturbed function $f_\tau\in\gt$
converging to $f$ as $\tau$ goes to $1$ unrestrictedly and such
that the solution of (\ref{aux1}) converges to the original
solution of equation (\ref{hgen}) uniformly on compact subsets of
$\Delta$.}

Geometrically, an affirmative answer to this question would give
us a constructive method for approximation of spirallike
(starlike) functions with respect to a boundary point by
spirallike (starlike) functions with respect to interior points.

We solve this problem as follows. Given $\tau$, we find a
transformation $\Phi_\tau:\Hol(\Delta,\C)\mapsto\Hol(\Delta,\C)$
which takes $h\in\spi[1]$ to $h_\tau\in\spi[\tau],\ \tau\in\Delta$
(i.e., $h_\tau=\Phi_\tau(h)$), and such that $\Phi_\tau(h)$ tends
to $h$ when $\tau$ tends to $1$ (see Theorem \ref{t3.4}).

To do this, we need first to consider some approximation and
interpolation problems for the classes $\G$ and $\Pp$ which are of
independent interest.

\vspace{2mm}

\subsection{ An approximation problem for generators}
As already mentioned, if a solution $h$ of equation (\ref{hgen})
is univalent, then $\mu$ must lie in the region
$\Omega=\biggl\{w\not=0:\ |w-\beta|\le\beta,\ \mbox{or}\
|w+\beta|\le\beta\biggr\}$, where $\beta=f'(\tau)>0$.

Actually, the instability phenomenon we have seen in Example~1
above follows from the fact that $\mu_\tau=f_\tau'(\tau)$ does not
necessarily converge to $\mu$ as $\tau\to1$ even if $\mu=f'(1)$.
Therefore, one can pose the following problem on controlled
approximation for functions in the class $\gt,\ \tau\in\pl\Delta$.

{\bf Let $f\in\gup$ with $f'(1)=\beta>0$. For $\tau\in\Delta$ and
given $\mu\in\Omega_+$, find a net $\{f_\tau\},\ f_\tau\in\gt,$
converging locally uniformly to $f$ as $\tau\to1$ unrestrictedly
and such that $\mu_\tau=f'_\tau(\tau)$ converges to $\mu$.}

We show how to solve this problem in Section 3.

\vspace{2mm}

\subsection{An approximation problem for semigroups}
For each $f\in\gt,\ \tau\in\overline\Delta,$ the generated
semigroup $\Ss=\left\{F_t\right\}_{t\ge0}\in\bt$ converges to the
point $\tau$ as $t\to\infty$. Moreover, the rate of convergence of
the semigroup can be estimated by the derivative $f'(\tau)$ of the
generator $f$ at the point $\tau$ (for the boundary case
$\tau\in\pl\Delta$, this fact follows by Theorem~\ref{ThD}; for
the general case $\tau\in\overline\Delta$, we refer to
\cite{E-R-S2002}, see also \cite{SD}).

Note that the number $\mu$ in Section 2.2 is not necessarily real,
while the angular derivative $f'(1)$ is a nonnegative real number.
So, in light of the instability phenomena mentioned above, the
following question seems to be natural.

{\bf Let $f\in\gup$ generate the semigroup
$\Ss=\left\{F_t\right\}_{t\ge0}\in\bu$. Suppose that a net
$\{f_\tau\},\ f_\tau\in\gt,\ \tau\in\Delta,$ converges to $f$
locally uniformly as $\tau\to1$. Do the corresponding elements
$F_{t,\tau},\ t\ge0,$ of the semigroups $\mathcal{S}_\tau$
generated by $f_\tau$ converge to $F_t$ for all $t\ge0$ as $\tau$
tends to $1$?}

Here the difficulty is that the value $F_{t,\tau}(z)$ must lie in
a neighborhood of $\tau\in\Delta$ for any $z\in\Delta$ and for $t$
sufficiently large, while $F_t(z)$ is close to $1$.

Nevertheless, we shall show that for each $r\in(0,1)$ and $T>0$,
the net $\left\{F_{t,\tau}(z)\right\}$ converges to $F_t(z)$
uniformly on the set $[0,T]\times(r\Delta)$.

\vspace{2mm}

\subsection{An approximation problem for functions with positive real part}
By the Berkson--Porta representation (see Theorem~\ref{ThC}), each
cone $\gt$, ${\tau\in\overline\Delta},$ can be parameterized by
elements of the cone $\Pp$ of functions with positive real part.
It turns out that the problems in Sections 2.1--2.3 lead us to the
following interpolation question about the approximation of
functions of class $\Pp$. Recall that a holomorphic function
$f\in\Hol(\Delta,\C)$ is said to be conformal at a boundary point
$\tau\in\partial\Delta$ if its angular derivative $f'(\tau)$
exists finitely and $f'(\tau)\not=0$.

{\bf Let $q\in\Pp$ be conformal at the boundary point $1$. For
$\tau\in\Delta$ and given
$\varphi\in\left(-\frac\pi2,\frac\pi2\right)$, does there exist a
net $\{q_\tau\}\subset\Pp$ converging locally uniformly to $q$ as
$\tau\to1$ unrestrictedly and such that}
\[
\arg q_\tau(\tau)=\varphi\quad\mbox{\it for all }\tau\in\Delta.
\]

A solution of this problem is given in Section 3.

\vspace{3mm}

\section{Main results}
\setcounter{equation}{0}

The well-known Riesz--Herglotz formula
\begin{equation}
p(z)=\oint\limits_{\partial\Delta}
\frac{1+z\overline{\zeta}}{1-z\overline{\zeta}}\,dm_{p}(\zeta)+i\Im
p(0)
\end{equation}
establishes a linear one-to-one correspondence between class $\Pp$
and the set of all positive measure functions $m(=m_p)$ on the
unit circle.

A consequence of this formula is that fact that for each
$\tau\in\pl\Delta$ the angular limit
\begin{equation}\label{delta}
\delta _{p}(\tau )=\angle \lim_{z\rightarrow 1}(1-z\overline{\tau }%
)p(z)=2m_{p}(\tau )
\end{equation}
exists and is a nonnegative real number.

\begin{defin}
The number $\delta=\delta _{p}(\tau )$ defined by (\ref{delta}) is
the {\sf charge} of the function $p\in\Pp$ at the boundary point
$\tau \in\pl\Delta$.
\end{defin}

This number plays a crucial role in our further considerations. In
particular, by using the Julia--Wolff--Carath\'eodory Theorem, one
can show that
\[
\delta_p(\tau)=2\inf_{z\in\Delta}\frac{1-|z|^2}{|1-z|^2}\,\Re p(z)
\]
(cf. \cite{SD1}).

We denote by $\pup$ the subclass of $\Pp$ consisting of functions
with positive charges at $\tau =1$, i.e., $p\in\pup$ if and only
if $\delta_{p}(1)>0$. Thus, $p\in\pup$ if and only if the function
$q\in\Pp$ defined by $q=\frac{1}{p}$ is conformal at $\tau =1$
with $q(1)=0$. Moreover, in this case, $q'(1)$ is a negative real
number.

The following assertion, which gives a solution to the problem in
Section~2.4, is the key for our further considerations.

\begin{theorem}
Let $q\in\Pp$ be conformal at $\tau =1$ with $q(1)=0$. Then for
all $\tau\in\Delta$ and each $\gamma\in\C$ such that
\begin{equation}\label{4}
\Re\gamma\ge\frac{\alpha}{2},\quad\mbox{where }\
\alpha=-q'(1)\left(>0\right),
\end{equation}
there exist functions
$\left\{q_\tau\right\}_{\tau\in\Delta}\subset\Pp$ converging to
$q$ uniformly on compact subsets of $\Delta$ when $\tau$ tends to
$1$ unrestrictedly and such that
\begin{equation}
q_\tau(\tau)=\gamma (1-\left|\tau\right|^{2})\rightarrow 0\ \mbox{
as }\ \tau\rightarrow 1 .
\end{equation}
In particular, if $\gamma$ is real and
$\dst\gamma\ge\frac{\alpha}{2}\,$, the values $q_\tau(\tau)$ are
real numbers.
\end{theorem}

\pr For $\gamma\in\C$, consider the function
\begin{equation}\label{6}
r(z)=\frac{zq(z)+\gamma z^{2}-\bar\gamma-2iz\Im\gamma}
{-(1-z)^{2}}.
\end{equation}
Clearly, $r\in\Hol(\Delta,\C)$ and $r(0)=\bar\gamma$.

We claim that $r\in\Pp$ if and only if $\gamma$ satisfies
inequality (\ref{4}).

Indeed, consider the function
\begin{equation}\label{7}
f(z)=zq(z)+\gamma z^2-\bar\gamma-2iz\Im\gamma,
\end{equation}
which is the numerator of (\ref{6}).

It follows by Theorem~\ref{ThC} that $f\in\G$. In addition,
\[
f(1)\left(=\angle\lim_{z\to1}
f(z)\right)=\gamma-\bar\gamma-2i\Im\gamma=0
\]
and
\bep
&&f'(1)\left(=\angle\lim_{z\to1}\frac{f(z)}{z-1}\right)\vspace{3mm}\\
&&=q'(1)+\angle\lim_{z\to1}\frac{\gamma
z^{2}-\bar\gamma-2iz\Im\gamma}{z-1} =q'(1)+2\Re\gamma .
\eep
Thus,
\begin{equation}\label{9}
f'(1)=2\Re\gamma -\alpha\geq 0
\end{equation}
if and only if condition (\ref{4}) holds.

At the same time, it follows from the infinitesimal version of the
Julia--Wolff--Carath\'{e}odory Theorem (Theorem~\ref{ThD}) that a
generator $f\in\G$ belongs to the subclass $\gu$ if and only if
$f(1)=0$ and $f'(1)\ge0$. On the other hand, by the uniqueness of
the Berkson--Porta representation (see Theorem~\ref{ThC}) of the
class $\G$, the inequality in (\ref{9}) is equivalent to
\[
\Re\frac{f(z)}{-(1-z)^2}\ge0.
\]
Comparing (\ref{9}) with (\ref{6}) and (\ref{7}) proves our claim.

Now for $\tau\in\Delta$, define the function $q_\tau$ on $\Delta$
by
\begin{equation}\label{qtau}
q_\tau(z)=\frac{1}{z} \left[\left(z-\tau\right)
\left(1-z\bar\tau\right) r(z) + \bar\gamma\tau-\gamma\bar\tau z^2
+ 2iz\Im\gamma\right]
\end{equation}
or
\begin{equation}\label{12}
q_\tau(z)=\frac{1}{z}g_\tau(z),
\end{equation}
where
\begin{equation}\label{13}
g_\tau(z)=\left(z-\tau\right) \left(1-z\bar\tau\right) r(z)
+\bar\gamma\tau-\gamma\bar\tau z^2+2iz\Im\gamma.
\end{equation}
Since $r(0)=\bar{\gamma}$, we have $g_\tau(0)=0$. Hence, $q_\tau$
is holomorphic in $\Delta$.

Moreover, since $\G$ is a real cone, we have by Theorem~\ref{ThC}
(ii)--(iii) that $g_\tau\in\G$ for all $\tau\in\Delta$.
Consequently, $g_\tau\in\G[0]$.

Again by Theorem~\ref{ThC} (ii), $\dst\Re
q_\tau(z)\left(=\Re\frac{g_\tau(z)}z\right)\geq 0$ for
$z\in\Delta$.

In addition, since
\begin{equation}\label{14}
q(z)=\frac{1}{z}\, g(z),
\end{equation}
where
\begin{equation}
g(z)=-\left(1-z\right)^2 r(z) - \gamma
z^2+\bar\gamma+2iz\Im\gamma,
\end{equation}
we have by (\ref{12}) and (\ref{13}) that $q_\tau(z)$ converges to
$q(z)$ as $\tau$ tends to $1$ unrestrictedly.

Finally, direct calculation shows that
$q_\tau(\tau)=\gamma\left(1-|\tau|^2\right)$.

\vspace{2mm}

\noindent{\bf Example 2.} Consider the function $q(z)=1-z$, which
obviously has positive real part. Then the function $r$ defined by
(\ref{6}) has the form
\[
r(z)=\frac1{1-z}\,.
\]
Substituting this function into formula (\ref{qtau}), we find the
approximating functions
\[
q_\tau(z)=\bar\tau(1-z)+\frac{|1-\tau|^2}{1-z}\,.
\]
We see that $q_\tau\to q$ as $\tau$ tends to $1$ unrestrictedly.

(i) Choosing the sequence of real numbers
$\tau_n^{(1)}=1-\frac1n$, we get the approximating sequence
$q_n^{(1)}$ of functions with positive real part. In Figure~1, we
see the images of the unit circle $\pl\Delta$ under the
approximating functions $q_1^{(1)},\ q_2^{(1)}$ and $q_4^{(1)}$ as
well as the image of $q$. It is worth noting that
$q_1^{(1)}(\Delta)=\{w:\, \Re w>\frac12\}$ and that the images
$q_n^{(1)}(\Delta)$ increase in the sense that
$q_n^{(1)}(\Delta)\subset q_{n+1}^{(1)}(\Delta)$. As $n$ tends to
infinity, these images tend to the right half-plane, while the
original function $q$ is bounded. (Note that the Carath\'eodory
Kernel Theorem is not applicable here since the functions $q_n$
are not univalent.)
\begin{figure}\centering 
    \includegraphics[angle=270,width=6cm,totalheight=6cm]{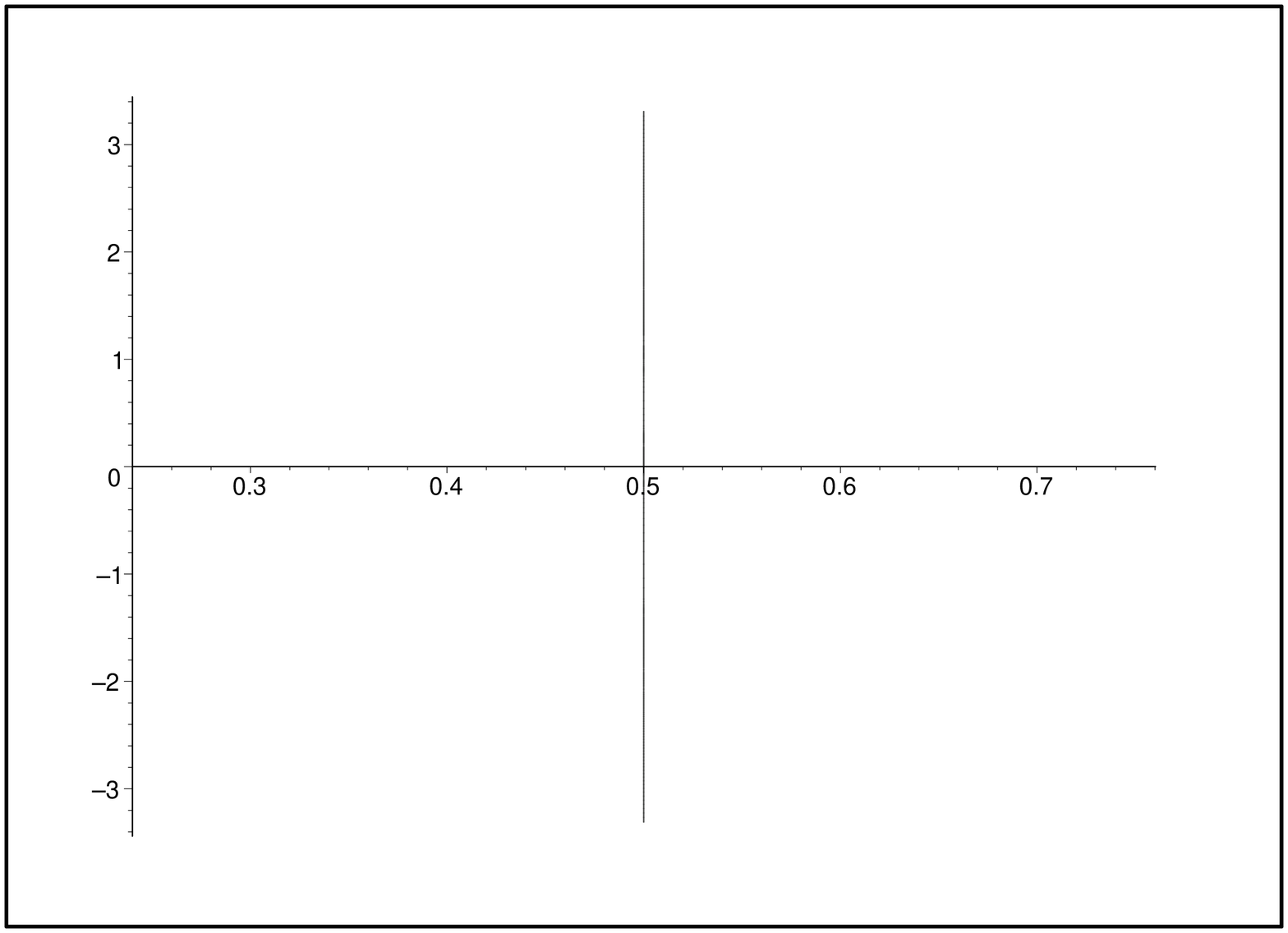}
  \hspace{3mm}
    \includegraphics[angle=270,width=6cm,totalheight=6cm]{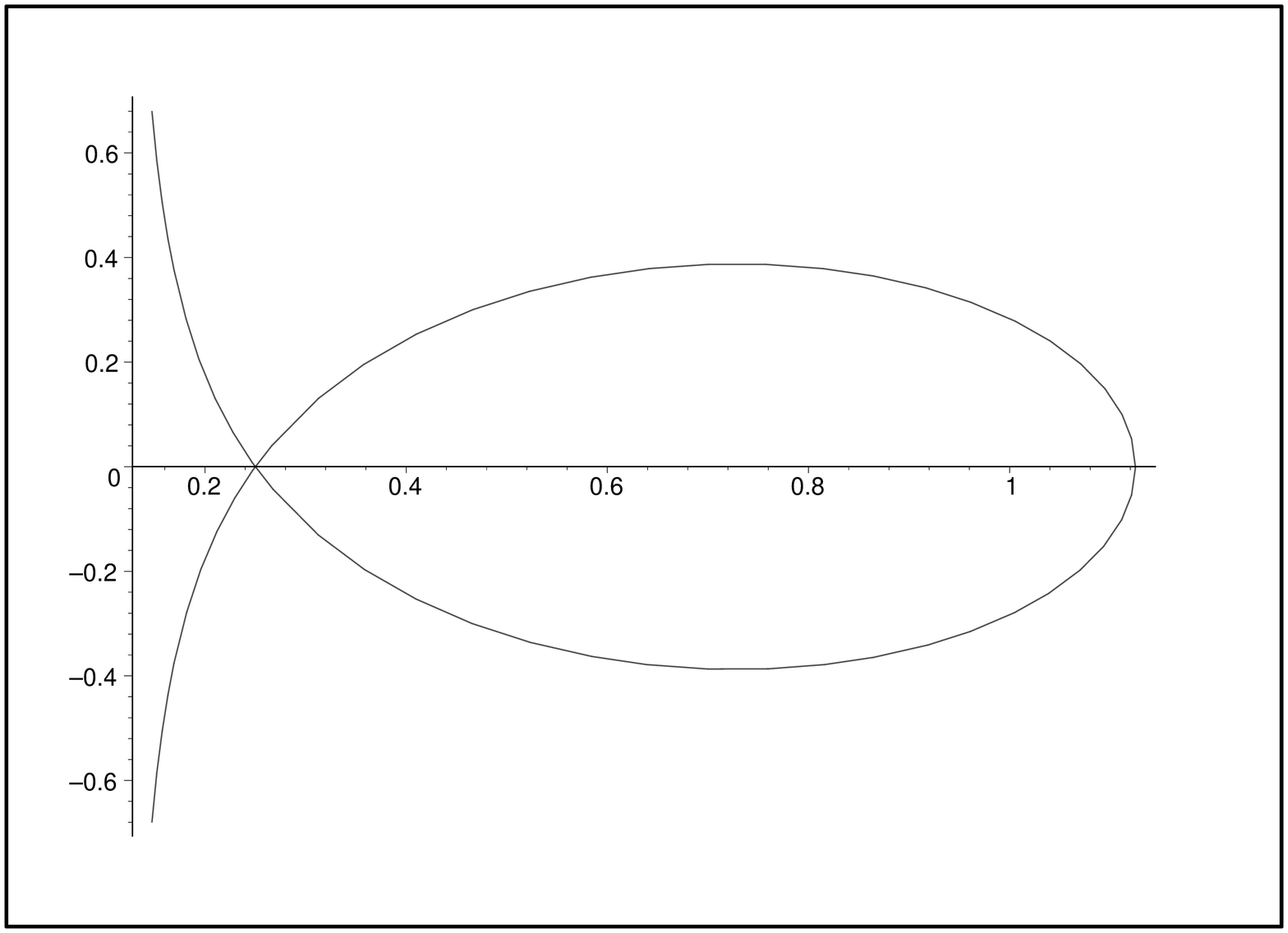}
\\ \vspace{3mm}
    \includegraphics[angle=270,width=6cm,totalheight=6cm]{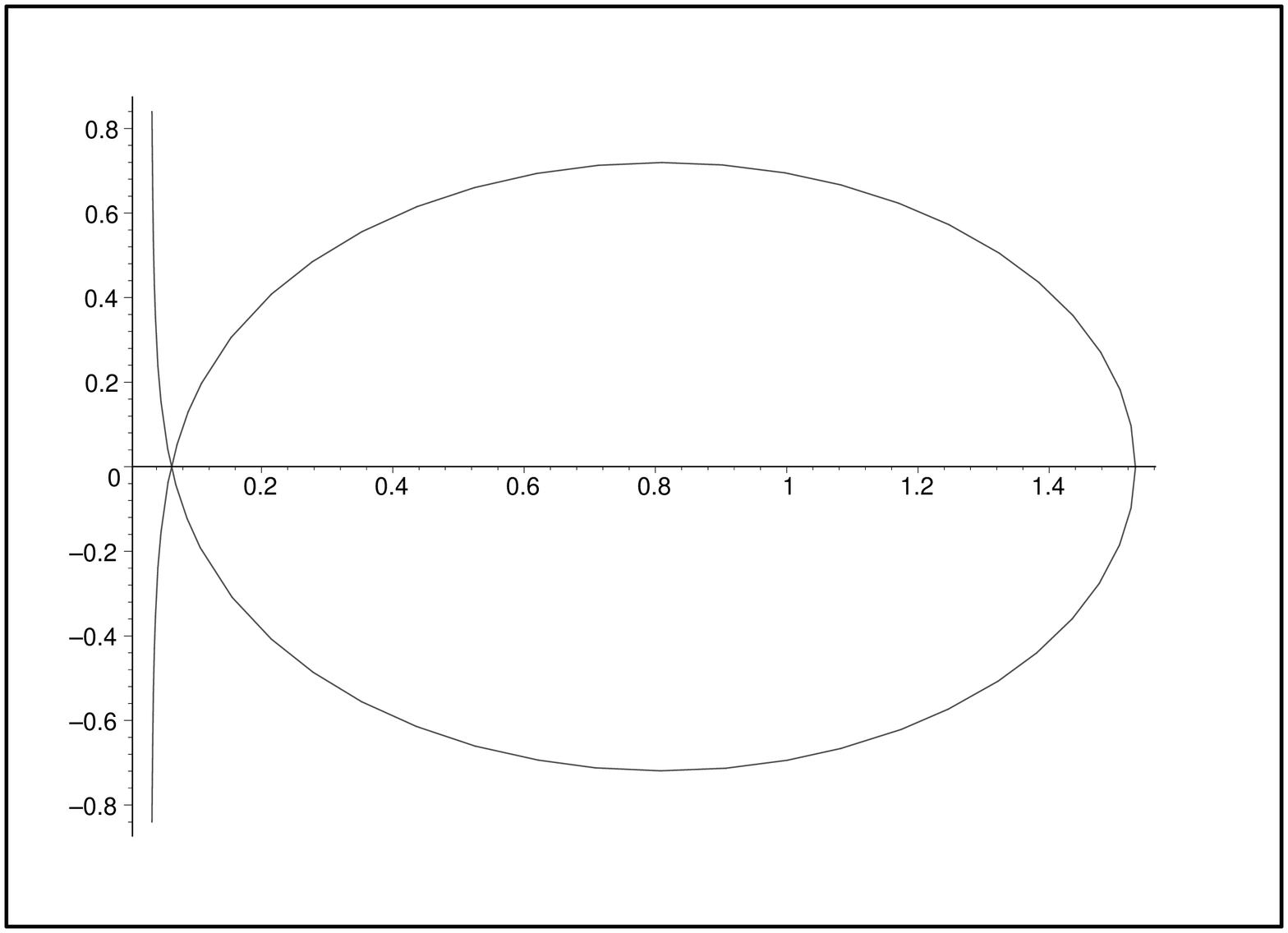}
\hspace{7mm}
    \includegraphics[angle=270,width=6cm,totalheight=6cm]{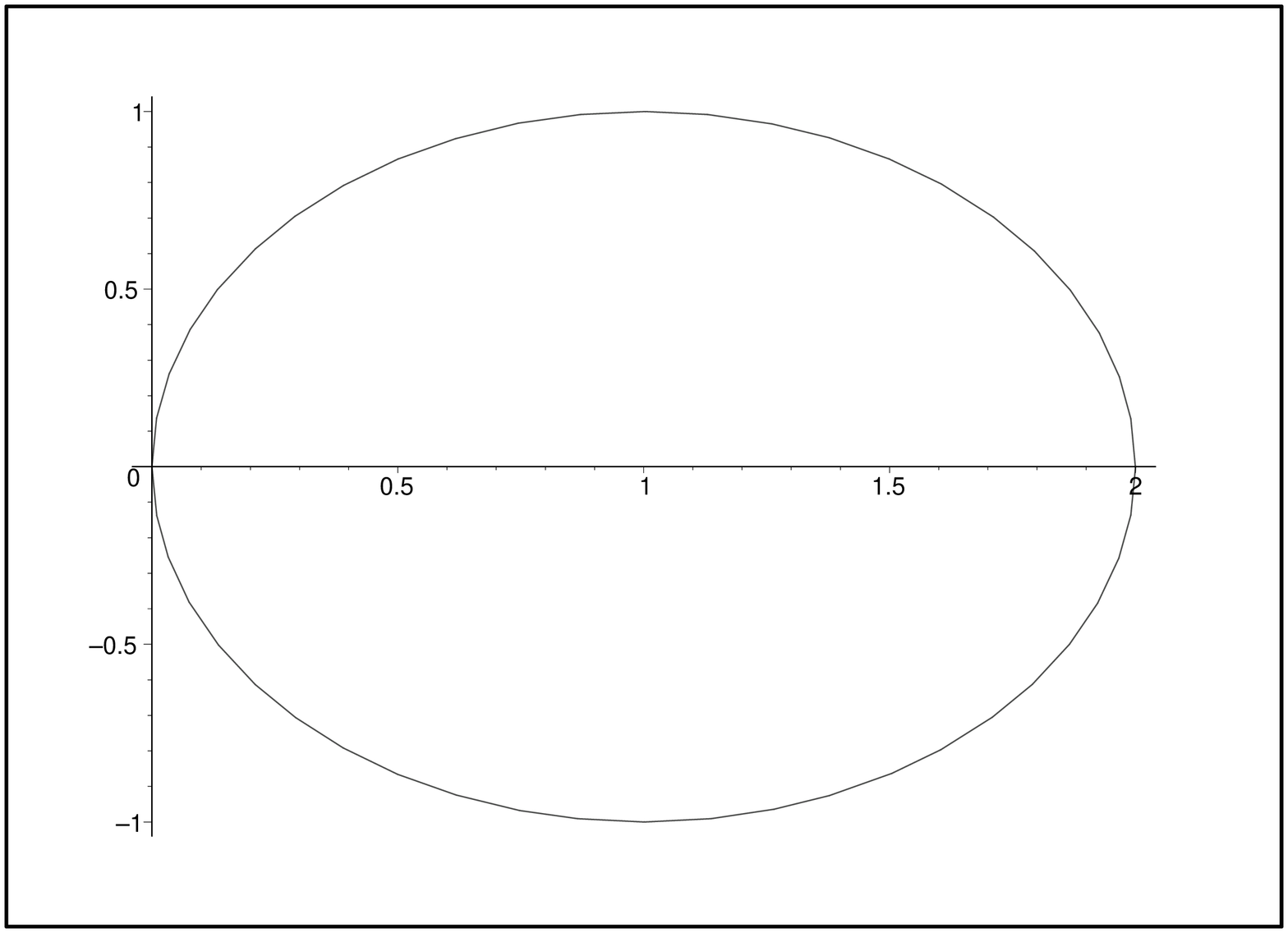}
    \caption{Example 2 (i), the images of  $q_1^{(1)},\ q_2^{(1)},\ q_4^{(1)}$ and of $q$.}
\end{figure}

(ii) Choosing another sequence $\tau_n^{(2)}$ converging to $1$,
say, $\tau_n^{(2)}=1-\frac{3(1-i)}{n}$, we get a different
approximating sequence $q_n^{(2)}$ of positive real part
functions. In Figure~2, we see the images of the unit circle
$\pl\Delta$ under the functions $q_4^{(2)},\ q_6^{(2)}$ and
$q_{12}^{(2)}$. Once again, in spite of the boundedness of $q$,
for $n$ large enough the image $q_n^{(2)}(\Delta)$ almost covers
the right half-plane.
\begin{figure}\centering 
    \includegraphics[angle=270,width=6cm,totalheight=6cm]{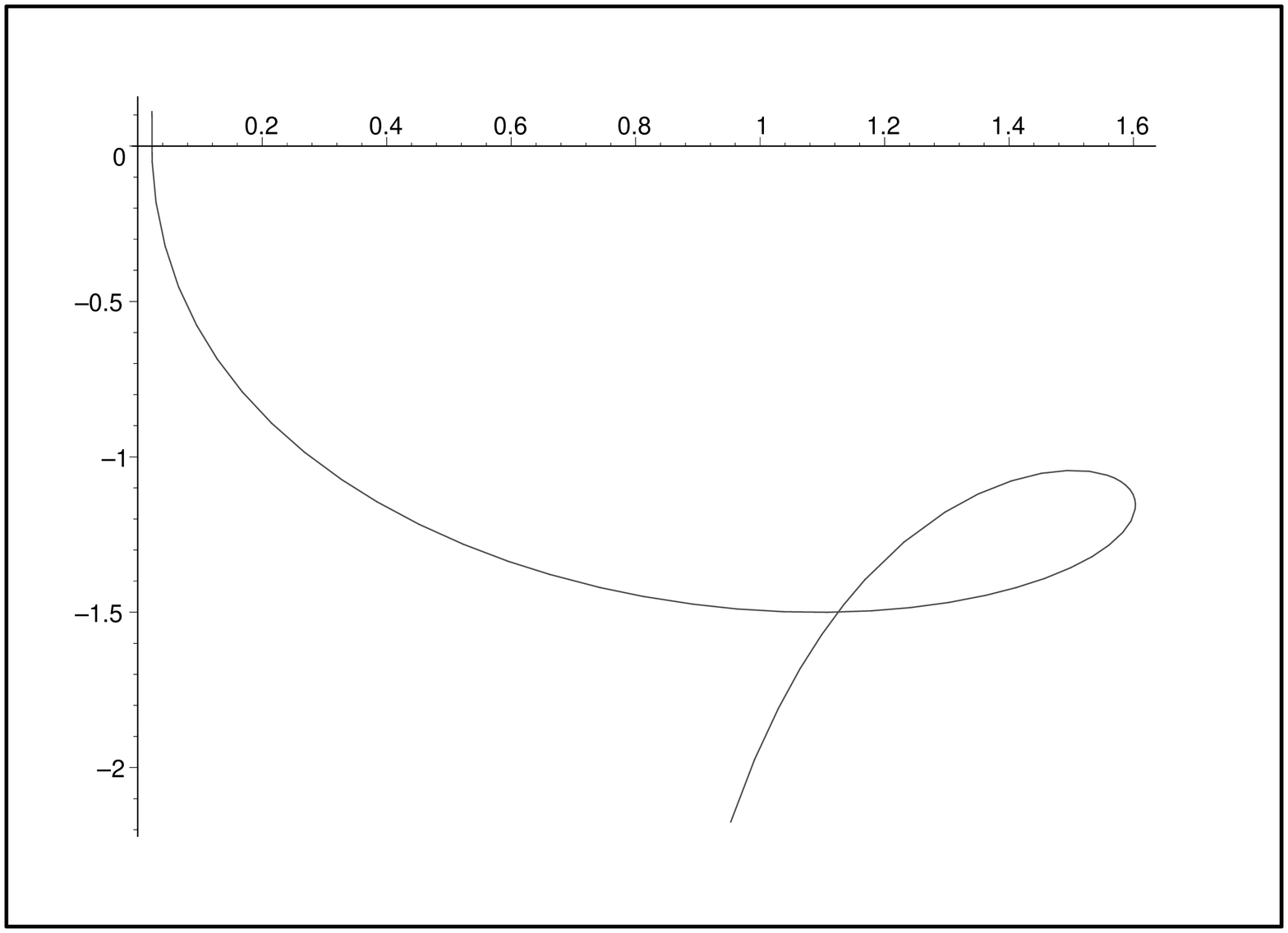}
  \hspace{3mm}
    \includegraphics[angle=270,width=6cm,totalheight=6cm]{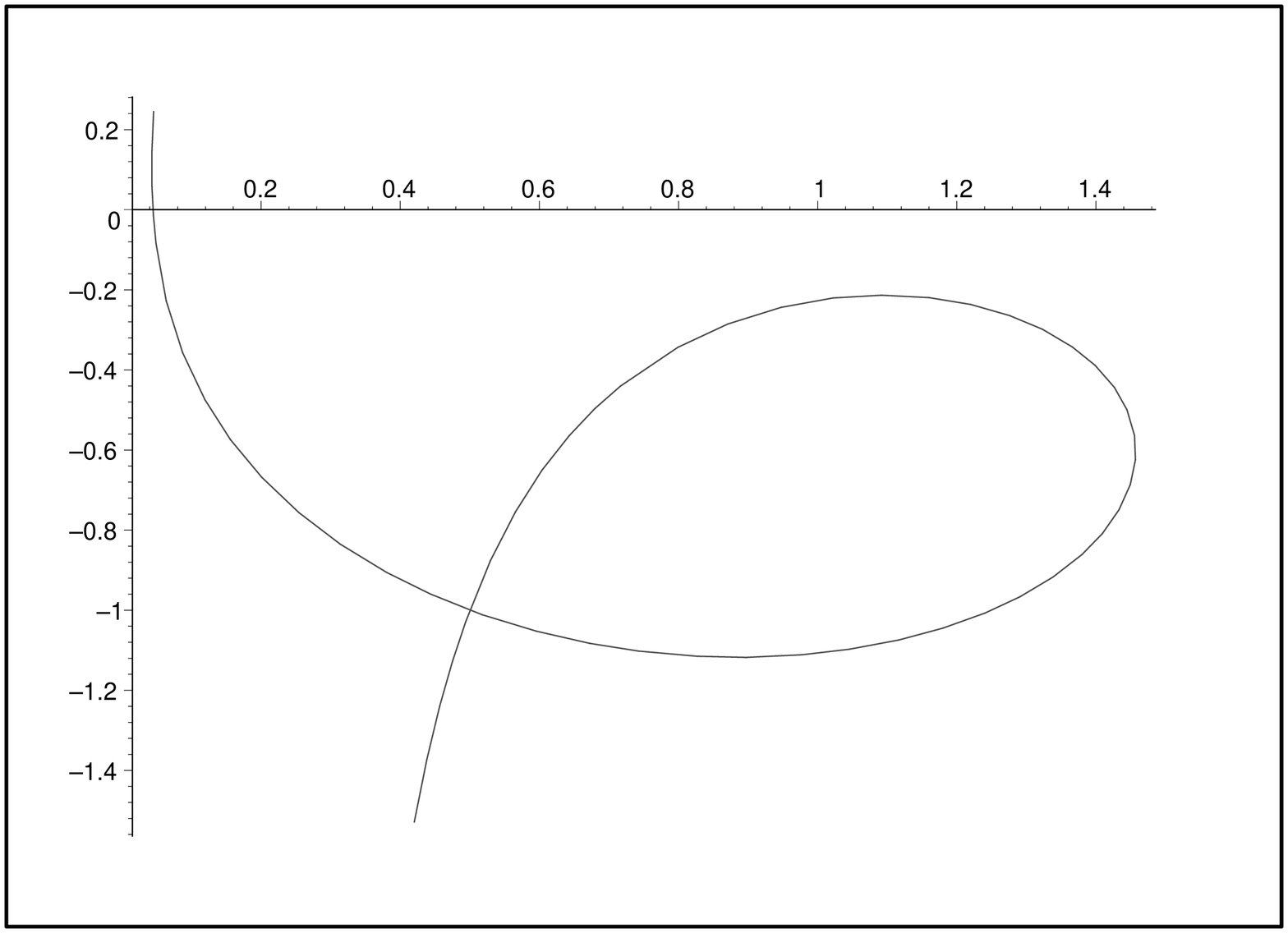}
\\ \vspace{3mm}
    \includegraphics[angle=270,width=6cm,totalheight=6cm]{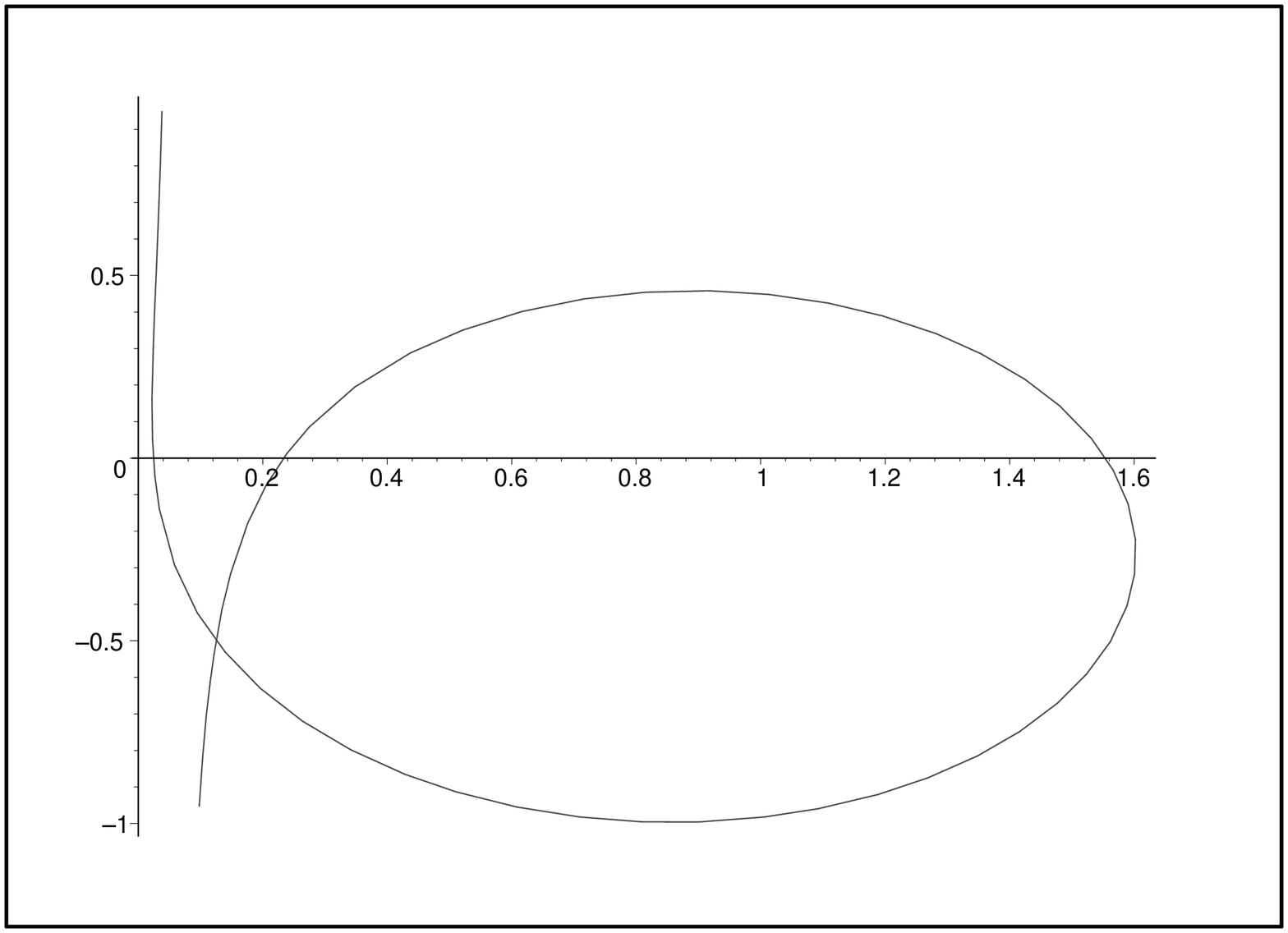}
    \hspace{7mm}
    \includegraphics[angle=270,width=6cm,totalheight=6cm]{f21.eps}
    \caption{Example 2 (ii), the images of  $q_4^{(2)},\ q_6^{(2)},\ q_{12}^{(2)}$.}
\end{figure}

\begin{corol}
Let $p\in\pup$ with the charge $\delta_p(1)=\beta >0$. Then for
all $\tau\in\Delta$ and each $\mu\in\C$ such that
\begin{equation}\label{3.14}
| \mu -\beta|\le \beta ,\ \mu \neq 0,
\end{equation}
there exist functions
$\left\{p_\tau\right\}_{\tau\in\Delta}\subset\Pp$ which converge
to $p$ uniformly on compact subsets of $\Delta$ as $\tau$ tends to
$1$ and such that
\begin{equation}\label{17}
p_\tau(\tau)=\frac{\mu }{1-|\tau|^2}\rightarrow \infty .
\end{equation}%
In particular, if $\mu$ is real and
\begin{equation}
0<\mu \leq 2\beta,
\end{equation}
the values $p_\tau(\tau)$ are real numbers.
\end{corol}

\pr Setting
\[
q(z)=\frac{1}{p(z)}\, ,
\]
we have $q\in\Pp$ with $q(1)=0$ and
\[
q'(1)=\angle\lim_{z\to 1}\frac{q(z)}{z-1}= \angle\lim_{z\to
1}\frac{1}{(z-1)p(z)}= -\frac{1}{\beta}\left(:=-\alpha \right) .
\]
Using Theorem 3.1 for each $\gamma\in\C$ such that
\begin{equation}\label{19}
2\Re\gamma\ge\frac{1}{\beta}
\end{equation}
and all $\tau\in\Delta$, one can find $q_\tau\in\Pp$ converging to
$q$ as $\tau\to1$ and such that
\[
q_\tau(\tau)=\gamma\left(1-|\tau|^2\right).
\]
Setting now
\begin{equation}
p_\tau(z)=\frac{1}{q_\tau(z)}
\end{equation}
and $\gamma =\frac{1}{\mu}\,,$\ we obtain that (\ref{3.14}) is
equivalent to (\ref{19}) and $p_\tau$ converges locally uniformly
to $p$.

Finally,
\[
p_\tau(\tau)=\frac{1}{q_\tau(\tau)}=
\frac{1}{\gamma\left(1-|\tau|^2\right)}= \frac{\mu}{1-|\tau|^2}\,,
\]
and we are done.

\vspace{2mm}

It is now easy to prove the following theorem, which solves the
problem in Section 2.2.

\begin{theorem}
Let $f\in\G$ satisfy the conditions $f(1)=0$ and $f'(1)=\beta>0$,
i.e., $f\in\G^+[1]$. Then for each $\mu\in\C$ such that
\begin{equation}
|\mu -\beta|\le\beta,\ \mu \neq 0
\end{equation}
and for all $\tau\in\Delta$, there exist generators $\{f_\tau\}$
with $f_\tau\in\gt$ such that $f_\tau'(\tau)=\mu$ and $f_\tau$
converges to $f$ uniformly on compact subsets of $\Delta$ as
$\tau$ tends to $1$ unrestrictedly.

In particular, setting $\mu =\beta$ we have
$f_\tau'(\tau)=\beta\in(0,\infty)$ for all $\tau\in\Delta$.
\end{theorem}

\pr By the Berkson--Porta formula (see Theorem~\ref{ThC}), $f$
must be of the form
\begin{equation}\label{20}
f(z)= -\left( 1-z\right)^2p(z),
\end{equation}
where $p\in\pup$ with $\delta_{p}(1)=\beta>0.$

By Corollary 3.1, for each $\mu$ satisfying (\ref{3.14}) and
$\tau\in\Delta$, one can find $p_\tau\in\Pp$ with
$p_\tau(\tau)=\dst\frac\mu{1-|\tau|^2}$ and such that $p_\tau$
converges to $p$ as $\tau\to1$.

Now we construct $f_\tau\in\G$ by the Berkson--Porta
representation (\ref{bp})
\[
f_\tau(z)=\left(z-\tau\right) \left(1-z\bar\tau\right) p_\tau(z).
\]
Since $f_\tau'(\tau)=\left(1-|\tau|^2\right) p_\tau(\tau)$, we
obtain our assertion.

\vspace{2mm}

In some sense a converse assertion is also true.

\begin{theorem}
Let $f\in\gup$, i.e., $f(1)=0$ and $f'(1)=\beta>0$, and let
$f_{n}\in\G$ vanish at $\tau_n\in\Delta$. Suppose that $f_n$
converges to $f$ uniformly on compact subsets of $\Delta$. Then

(i) the sequence $\left\{\tau_n\right\}$ converges to $1$;

(ii) if the sequence $\mu_n:=f'_n(\tau_n)$ converges to
$\mu\neq0$, then $\mu$ must satisfy condition (\ref{3.14}), i.e.,
$|\mu-\beta|\le\beta$.
\end{theorem}

\pr (i) For each $n=1,2,\ldots,$ the function $f_n\in\G[\tau_n]$
has the form $f_n(z)=\left(z-\tau_n\right)\left(
1-z\overline{\tau_n}\right) p_n(z)$ with $\Re p_n(z)\ge0$
everywhere.

Since $\overline{\Delta}$ is a compact subset of $\C$ and
$\left\{p_n\right\}$ is a normal family on $\Delta$, one can
choose a subsequence $\left\{ n_{k}\right\}\subset\N$ such that
$\tau_{n_k}$ converges to $\tau\in\overline\Delta$ and $p_{n_k}$
converges either to $\widetilde{p}\in\Pp$ or
$\widetilde{p}=\infty$. In any case, the sequence
$\left\{f_{n_k}\right\}\subset\G$ converges to either
$f(z)=\left(z-\tau\right) \left(1-z\bar\tau\right)
\widetilde{p}(z)$ with $\Re\widetilde{p}(z)\ge0$ and
$\tau\in\overline\Delta$ or to infinity. Since the latter case is
impossible, the uniqueness of the Berkson--Porta representation
and (\ref{20}) imply that $\tau=1$ and $\widetilde{p}=p$.

(ii) Assume now that $\mu_n:=f_{n}'(\tau_n)$ converges to
$\mu\in\C,\ \mu\neq0$, and consider the differential equations
\begin{equation}
\mu_n h_n(z) = h_n'(z) f_n(z) ,
\end{equation}
normalized by the conditions $h_n(0)=1$. This initial value
problem has a unique solution which is a univalent function on
$\Delta $\ spirallike with respect to an interior point and
$h_n(\tau_n)=0\in h_{n}(\Delta)$.

Explicitly, $h_{n}$ can be written as
\begin{equation}\label{22}
h_n(z) = \exp \left[\mu_n\int_0^z\frac{dz}{f_n(z)}\right] .
\end{equation}
Now for each $r\in(0,1)$, one can find $n_{0}\in\N$ such that for
all\ $n>n_{0}$, the points $\tau_n$ do not belong to the closed
disk $\overline{\Delta_r}=\left\{ |z|\leq r<1\right\}$. In other
words, for all $n>n_0$ the functions $f_n$ do not vanish in this
disk. Then the functions $h_{n}$ defined by (\ref{22}) converge to
a function $h\in\Hol(\Delta_r,\C)$ uniformly on this disk. Since
$r$ is arbitrary, we have that actually $h\in\Hol(\Delta,\C)$ and
has the form
\begin{equation}\label{23}
h(z) =\exp\left[ \mu\int_0^z\frac{dz}{f(z)} \right] .
\end{equation}
In addition, it follows by Hurwitz's theorem that $h$ is either a
univalent function on $\Delta $ or a constant. The latter case is
impossible because of the equality
$\dst\frac{h'(z)}{h(z)}=\frac{\mu}{f(z)}\,,$ which follows by
(\ref{23}) and $\mu\neq 0$.

On the other hand, we already know that the initial value problem
\[
\beta \widetilde{h}(z) =\widetilde{h}'(z)
f(z),\quad h(0) =1
\]
has also the unique solution
\begin{equation}\label{24}
\widetilde{h}(z) =\exp \left[ \beta \int_{0}^{z}\frac{dz}{f(z)}
\right],
\end{equation}
which is a starlike function with respect to a boundary point
$\left(\widetilde{h}(1)=0\right)$.

Comparing (\ref{23}) and (\ref{24}), we obtain
\begin{equation}
h(z) = \left[ \widetilde{h}(z) \right] ^{\frac{\mu }{\beta }}.
\end{equation}
In addition, the Visser--Ostrowski condition
\[
\angle\lim_{z\to 1}\frac{( z-1)
\widetilde{h}'(z)}{\widetilde{h}(z)}=1
\]
(see \cite{E-R-S2001a}) implies that the smallest wedge which
contains $h(\Delta)$ is exactly of angle $\pi$. Since $h$ is
univalent, we get that $\dst\Re\frac{\beta }{\mu
}\ge\frac{1}{2}\,$, which is equivalent to~(\ref{3.14}).

\vspace{2mm}

As a consequence of Theorems 3.2 and 3.3, we obtain the following
result.

\begin{theorem}\label{t3.4}
Let $f\in\gup$ with $f'(1)=\beta$. Then the initial value problem
\begin{equation}\label{26}
\lambda h(z) =h'(z)f(z),\quad h(0) =1
\end{equation}
has a univalent solution $h\in\Hol(\Delta,\C)$ if and only if the
complex number $\lambda$ belongs to the set
$\Omega=\Omega_+\cup\Omega_-$, where
\begin{equation}
\Omega_\pm=\left\{\omega\in\C:\ |\omega\mp\beta|\le\beta,\
\omega\neq 0\right\} .
\end{equation}
Moreover,

\noindent$\bullet $ for each $\mu\in\Omega_+$ and each
$\tau\in\Delta$, there are $h_\tau\in\spi[\tau]$ which converge to
the function
\begin{equation}
h_{1}(z)=h^{\frac{\mu }{\lambda }}(z)
\end{equation}
as $\tau$ tends to $1$ unrestrictedly. In particular, if
$\lambda\in\Omega_+$, then choosing $\mu=\lambda$, we find
$h_\tau\in\spi[\tau]$ which converge to the original function $h$
as $\tau\to1$.

\noindent$\bullet$ for each $\mu\in\Omega_-$ and for each
$\tau\in\Delta$, there are meromorphic functions
$\widetilde{h}_\tau$ with a unique simple pole at $\tau$ and such
that $\widetilde{h}_\tau$ converges to the holomorphic function
\begin{equation}
\widetilde{h}(z)=h^{\frac{\mu }{\lambda }}(z)
\end{equation}
as $\tau$ tends to $1$ unrestrictedly.
\end{theorem}

\pr Take any $\mu\in\Omega_+$ and any $\tau\in\Delta$. By
Theorem~3.2, one can choose generators $f_\tau\in\gt$ such that
the functions $f_\tau$ converge uniformly on compact subsets to
$f$ and satisfy the conditions $f_\tau(\tau)=0$ and
$f_\tau'(\tau)=\mu$. Then, as in the proof of the second part of
Theorem~3.3, one shows that functions $h_\tau\in\spi[\tau]$
defined by
\begin{equation}
h_\tau=\exp \left[\mu\int_0^z \frac{dz}{f_\tau(z)}\right]
\end{equation}
converge to a univalent function
\begin{equation}\label{31}
h_1(z) =\exp\left[\mu\int_0^z\frac{dz}{f(z)}\right]
\end{equation}
which satisfies the equation
\begin{equation}\label{32}
\mu h_1(z)=h_1'(z) f(z).
\end{equation}
If now $\lambda\in\Omega_+$, then setting $\mu =\lambda$, we see
that $h=h_1$ must be univalent. If $\lambda\in\Omega_-$, then
setting $\mu =-\lambda$, and comparing differential equations
(\ref{26}) and (\ref{32}), we see that $h=h_1^{-1}$. Since
$h_1(z)\neq 0,\ z\in \Delta$, $h$ is a well-defined univalent
function on $\Delta$. In addition, it is clear that $h_{1}^{-1}$
is a locally uniform limit of meromorphic functions $h_\tau^{-1}$
with poles at $\tau$.

Conversely, assume that for some $\lambda\in\C,\ \lambda\neq 0$,\
equation (\ref{26}) has a univalent solution in $\Delta$. If
$\Re\lambda>0$, then setting $\mu=\beta$ in (\ref{31}), we see as
in the proof of Theorem 3.3 that the image of the function
$h_{1}(z)=h^{\frac{\beta}{\lambda}}(z)$ must lie in the wedge of
the angle $\pi$, which is the smallest one containing
$h_{1}(\Delta)$. Then, by a result of \cite{A-E-S}, we have
$\lambda\in\Omega_+$. If $\Re\lambda<0$, then the same
considerations show that $-\lambda\in\Omega_+$, and we are done.

\vspace{2mm}

\begin{theorem}[Perturbation formula]\label{pert}
Let $f\in\gup$ with $f'(1)=\beta$ and let
$\lambda\in\Omega_+=\left\{w\in\C:\ |w-\beta|\leq\beta,\
w\neq0\right\}$.

Assume that $h\in\Hol(\Delta,\C)$ is the solution of
equation~(\ref{26})
\[
\lambda h(z)=h'(z) f(z)
\]
normalized by the conditions $h(0)=1$ and $h(1)=0.$ Then for each
$\mu\in\Omega_+$ and for each $\tau\in\Delta$, the function
$h_\tau\in\Hol(\Delta,\C)$ defined by
\[
h_\tau(z)= \left[h(z)\right]^{\frac{\mu}{\lambda}}
\,\frac{(z-\tau)(1-z\bar\tau)^{\mu/\bar\mu}}
{-(1-z)^{1+\mu/\bar\mu}}
\]
is univalent on $\Delta$ and belongs to the class $\spi[\tau]$
with $h_\tau(\tau)=0$ and $h_\tau(0)=\tau .$

If, in particular, $\mu =\lambda$, then $h_\tau$ converges to $h$
whenever $\tau$ tends to $1.$ Thus $h$ is a univalent function on
$\Delta$ spirallike (starlike) with respect to a boundary point
with $h(1)=0.$
\end{theorem}

\pr Let $h\in\Hol(\Delta,\C)$ be a solution of the differential
equation (\ref{26}) with $f\in\gup$ defined by $\dst
f(z)=\frac{-(1-z)^2}{q(z)}$, where
\[
\Re q(z)>0 \quad\mbox{and}\quad
\angle\lim_{z\to1}\frac{1-z}{q(z)}=\beta.
\]
Then
\be\label{q}
\frac{h'(z)}{\lambda h(z)}=\frac{q(z)}{-(1-z)^2}\ .
\ee

As in the proof of Theorem 3.1, for $\gamma$ complex with
$2\Re\gamma\ge\frac1\beta$ we define $r\in\pup$ by
\be\label{r1}
r(z)=\frac{zq(z)+(\gamma z+\bar\gamma)(z-1)}{-(1-z)^2}\, .
\ee

Comparing (\ref{q}) and (\ref{r1}), we have
\be\label{r}
r(z)=\frac{zh'(z)}{\lambda h(z)}+\frac{\gamma\,
z+\bar\gamma}{1-z}\,.
\ee

Once again, for each point $\tau\in\Delta$, we define
\[
q_\tau(z)=\frac1z\left[(z-\tau)(1-z\bar\tau)r(z)+\bar\gamma\tau-\gamma\bar\tau
z^2+(\gamma-\bar\gamma)z \right].
\]

Then the differential equation
\be\label{htau1}
\mu h_\tau(z)=h_\tau'(z)(z-\tau)(1-z\bar\tau)\frac1{q_\tau(z)}
\ee
has a holomorphic solution $h_\tau$ if and only if
$\dst\mu=\frac1\gamma\in\Omega_+$. Moreover, $h_\tau\in\spi[\tau]$
by Theorems~1.4 and 1.2. Since $h_\tau(\tau)=0$, the function
$\dst g_\tau(z):=\frac{h_\tau(z)}{z-\tau}$ is well-defined. Now,
using (\ref{r}) and (\ref{htau1}), we calculate
\bep
\frac{g_\tau'(z)}{g_\tau(z)}=\mu\left[\frac{h'(z)}{\lambda h(z)}+
\frac{\gamma z+\bar\gamma}{z(1-z)} -
\frac{\bar\gamma}{z(1-z\bar\tau)} \right] \vspace{3mm} \\
=\mu\left[\frac{h'(z)}{\lambda h(z)}+ \frac{\gamma}{1-z}+
\frac{\bar\gamma}{1-z} + \frac{\bar\gamma(-\bar\tau)}{1-z\bar\tau}
\right] \vspace{3mm} \\
=\frac{\mu h'(z)}{\lambda h(z)}+ \frac1{1-z}+
\frac{\mu\bar\gamma}{1-z} +
\frac{\mu\bar\gamma(-\bar\tau)}{1-z\bar\tau} \,.
\eep

Then
\[
\dst h_\tau(z)=(z-\tau)g_\tau(z)=C\left(h(z)
\right)^{\frac\mu\lambda}\cdot
\frac{(z-\tau)(1-z\bar\tau)^{\frac\mu{\bar\mu}}}{(1-z)^{1+\frac\mu{\bar\mu}}}
\]
with some constant $C$. So, to satisfy the normalization
$h_\tau(0)=\tau$, we must set $C=-1$. The proof is complete.

\vspace{2mm}

\noindent{\bf Example 2.} Consider the starlike function with
respect to a boundary point
\[
h(z)=(1-z)^{0.8}.
\]
It satisfies equation~(\ref{26}) with $f(z)=z-1$ and
$\lambda=0.8$:
\[
0.8h(z)=h'(z)\cdot(z-1) .
\]
Setting $\mu=\lambda=0.8$, we see that by Theorem~\ref{pert}, $h$
can be approximated by the functions
\[
h_\tau(z)=\frac{(\tau-z)(1-z\bar\tau)}{(1-z)^{1.2}}.
\]

(i) Choosing, in particular, the sequence of real numbers
$\tau_n^{(1)}=1-\frac3n$, we get the approximating sequence
$h_n^{(1)}$ of starlike functions with respect to (different)
interior points. In Figure~3, we see the images of $h$ as well as
the images of the approximating functions $h_6^{(1)},\
h_{10}^{(1)}$ and $h_{30}^{(1)}$. Note that in this case, the
intersection $\bigcap_{n=0}^\infty h_n^{(1)}(\Delta)$ contains the
left half-plane, while the image $h(\Delta)$ lies in the right
half-plane.

\begin{figure}\centering 
    \includegraphics[angle=270,width=9.0cm,totalheight=9.0cm]{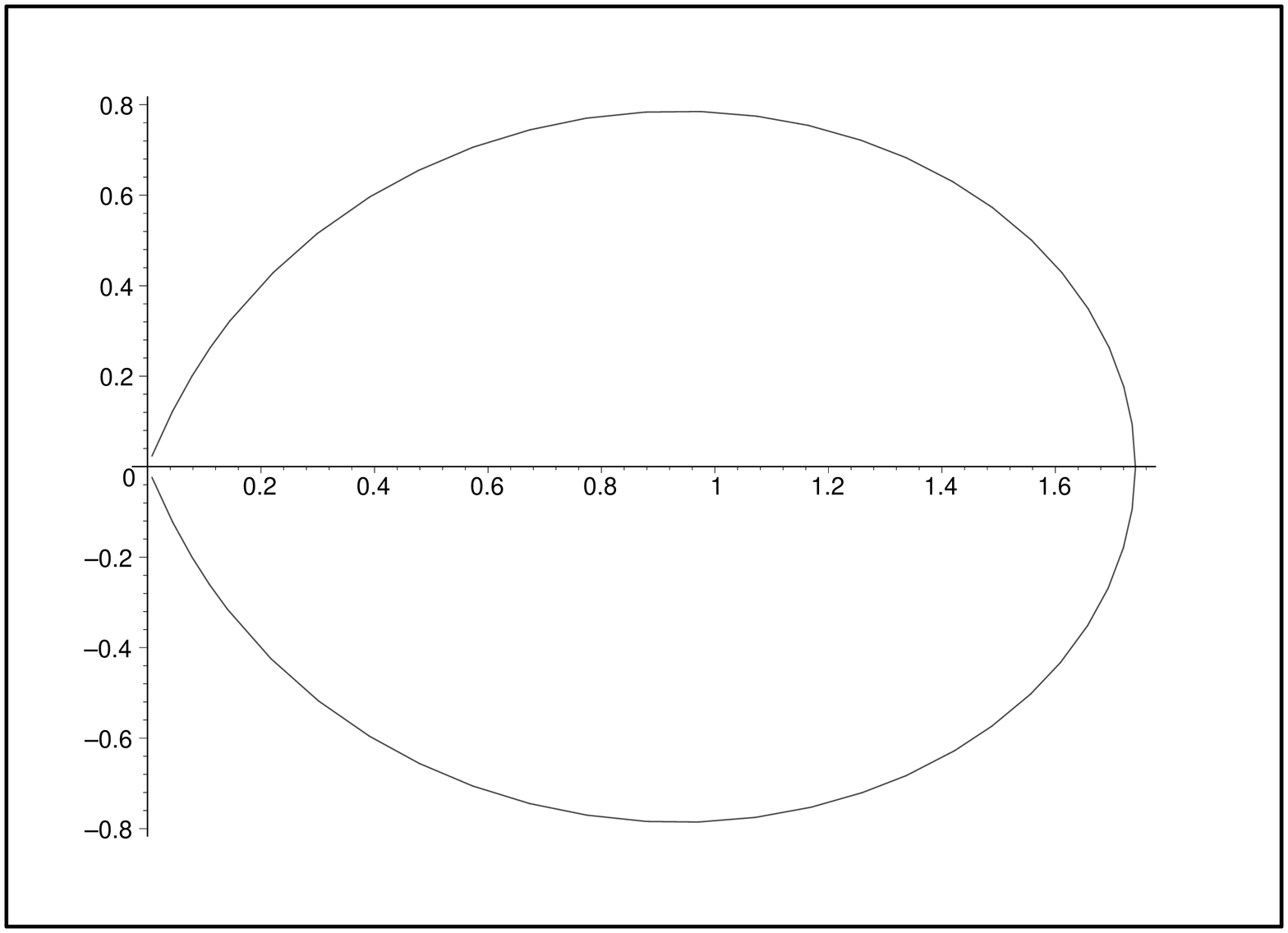}
  \hspace{3mm}
    \includegraphics[angle=270,width=3.8cm,totalheight=9cm]{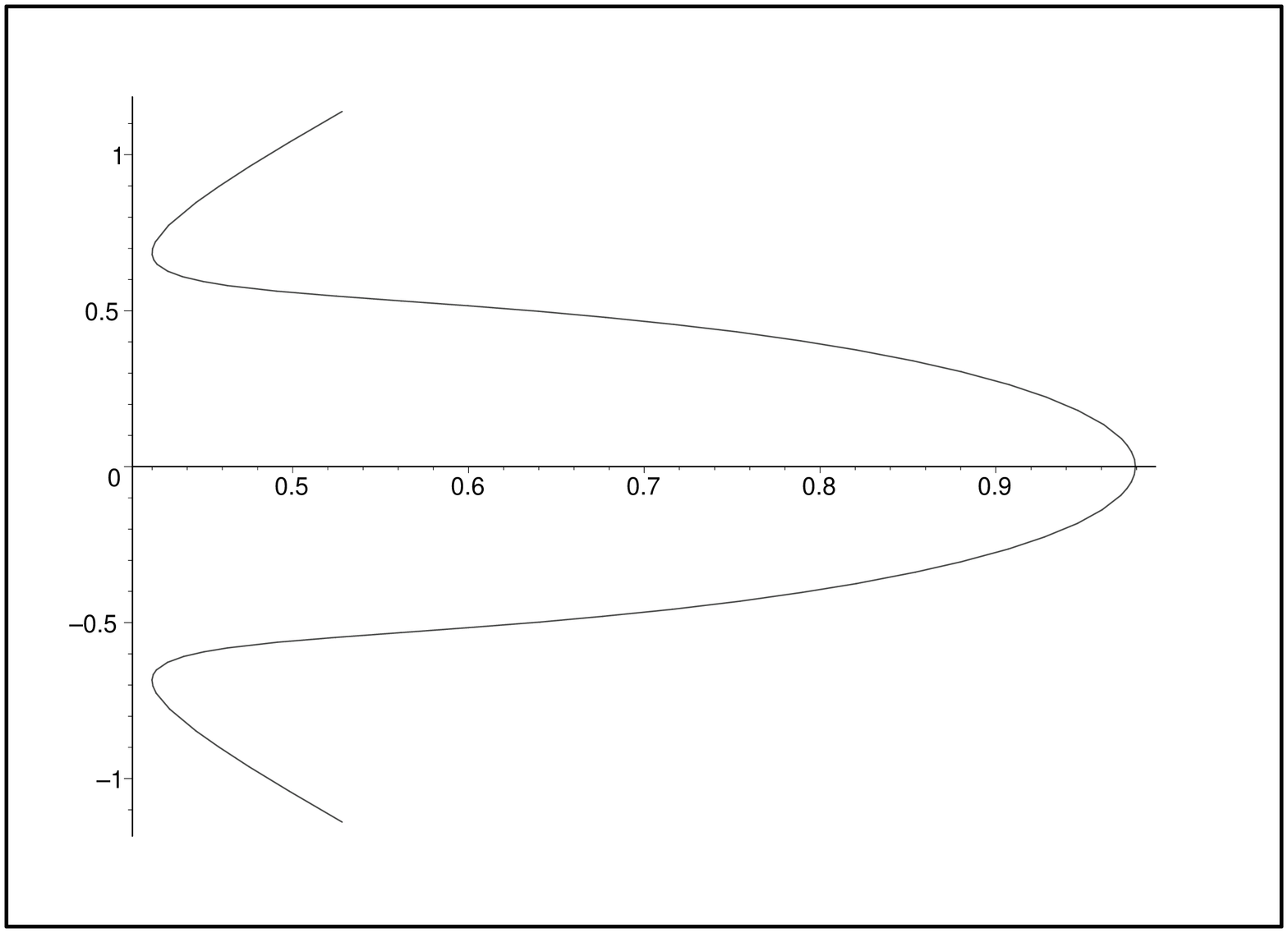}
\\ \vspace{3mm}
    \includegraphics[angle=270,width=4.7cm,totalheight=9.0cm]{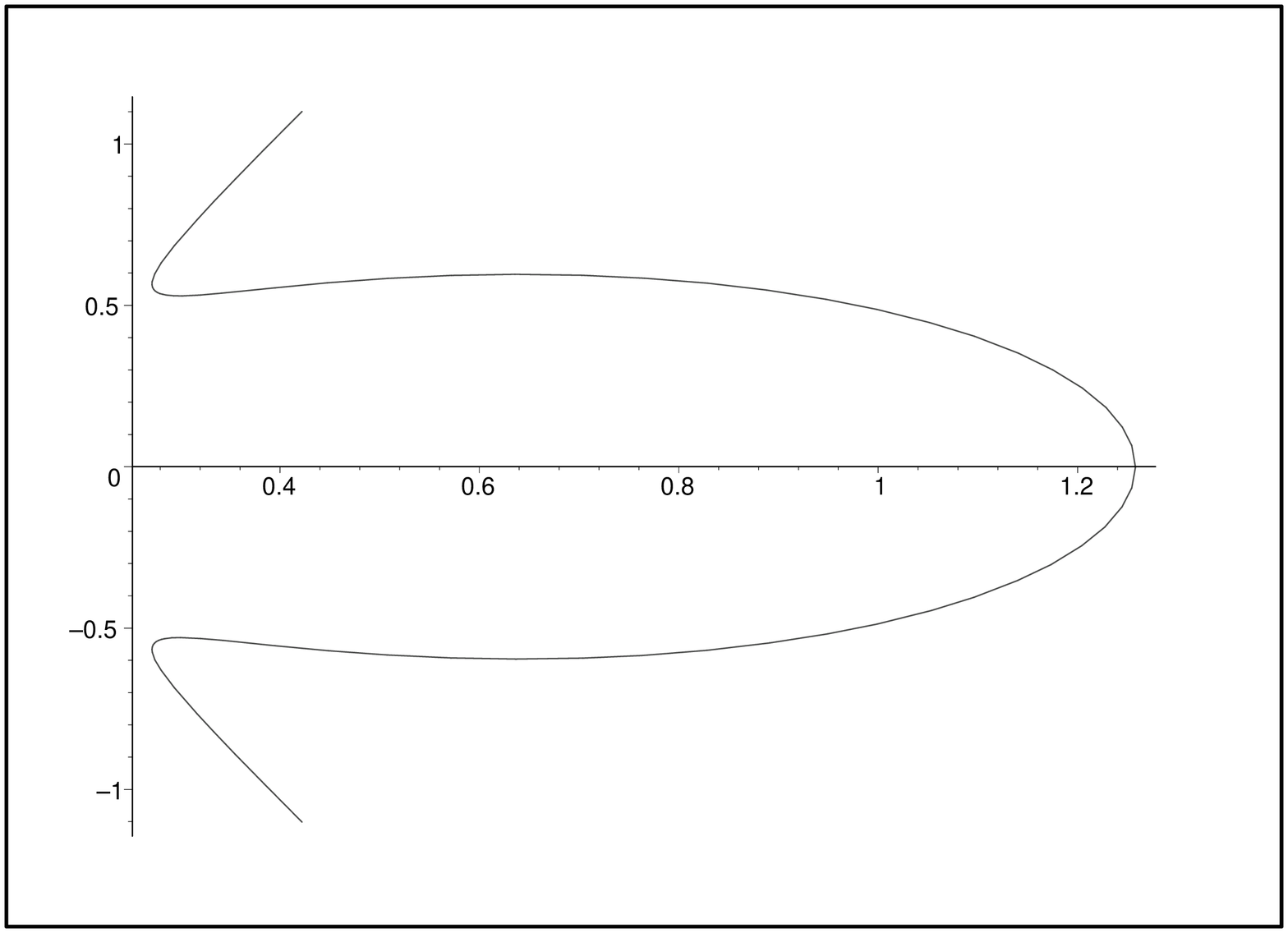}
\hspace{7mm}
    \includegraphics[angle=270,width=7cm,totalheight=9.0cm]{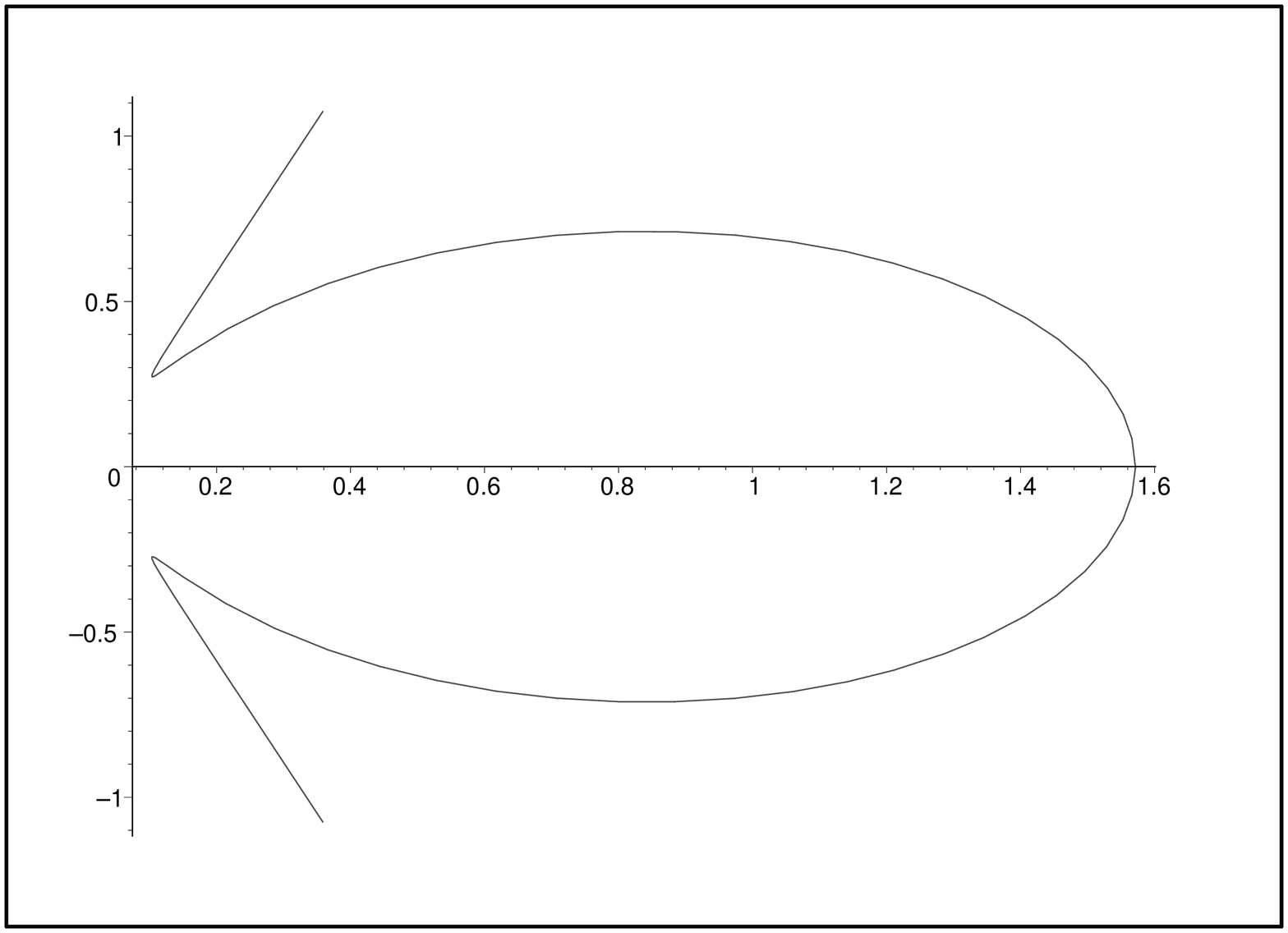}
    \caption{Example 3 (i), the images of $h, h_6^{(1)}, h_{10}^{(1)}$ and $h_{30}^{(1)}$.}
\end{figure}

(ii) Choosing another sequence $\tau_n^{(2)}$ converging to $1$,
say, ${\tau_n^{(2)}=1-(1-i)\frac3n}$, we get a different
approximating sequence $h_n^{(2)}$ of starlike functions with
respect to interior points. In Figure~4, we see the images of $h$
and of the approximating functions $h_6^{(2)},\ h_{10}^{(2)}$ and
$h_{30}^{(2)}$. Once again, all of the images $h_n^{(2)}(\Delta)$
contain the left half-plane.

\begin{figure}\centering 
    \includegraphics[angle=270,width=3.6cm,totalheight=9.0cm]{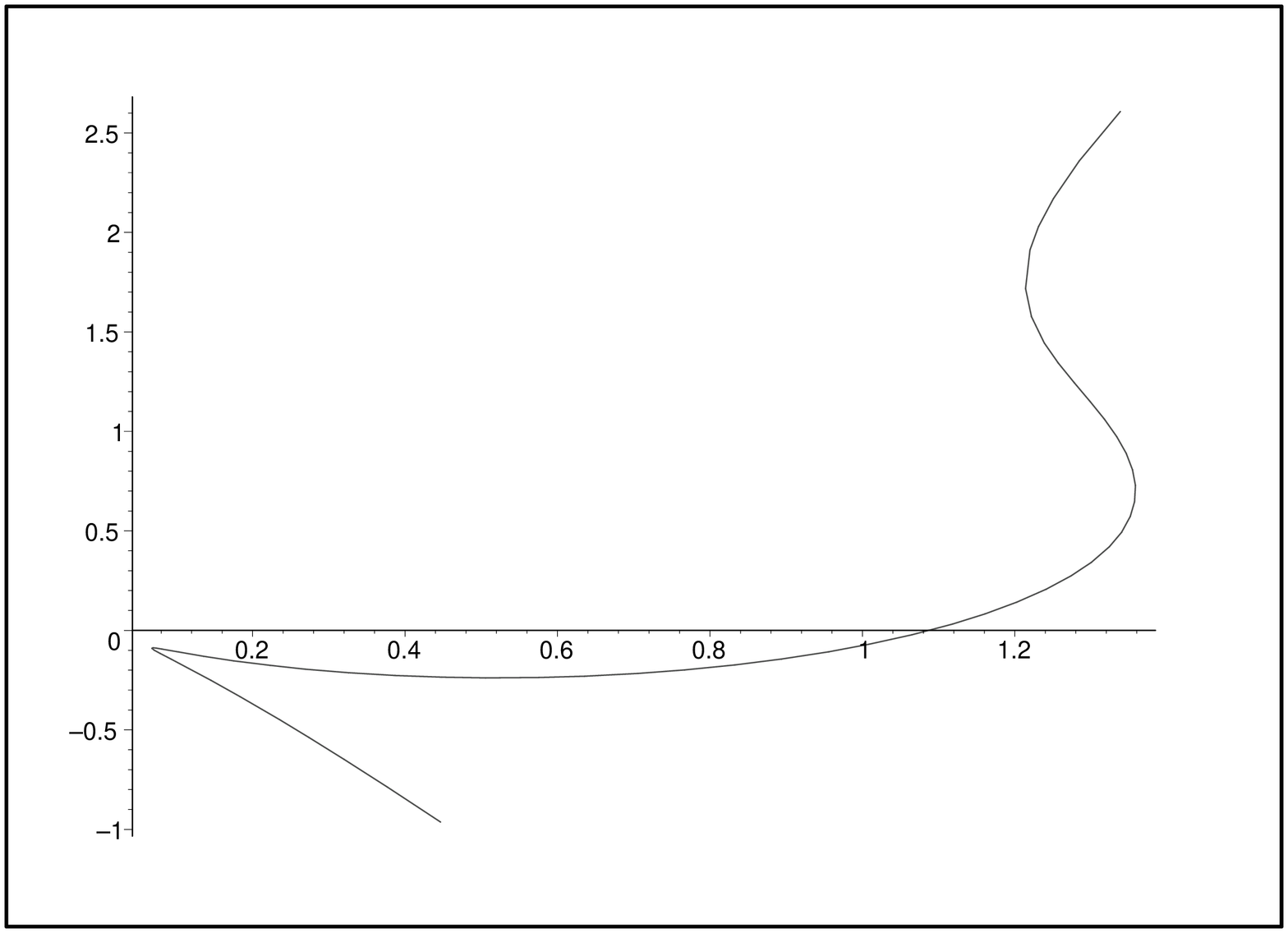}
\hspace{7mm}
    \includegraphics[angle=270,width=4.3cm,totalheight=9.0cm]{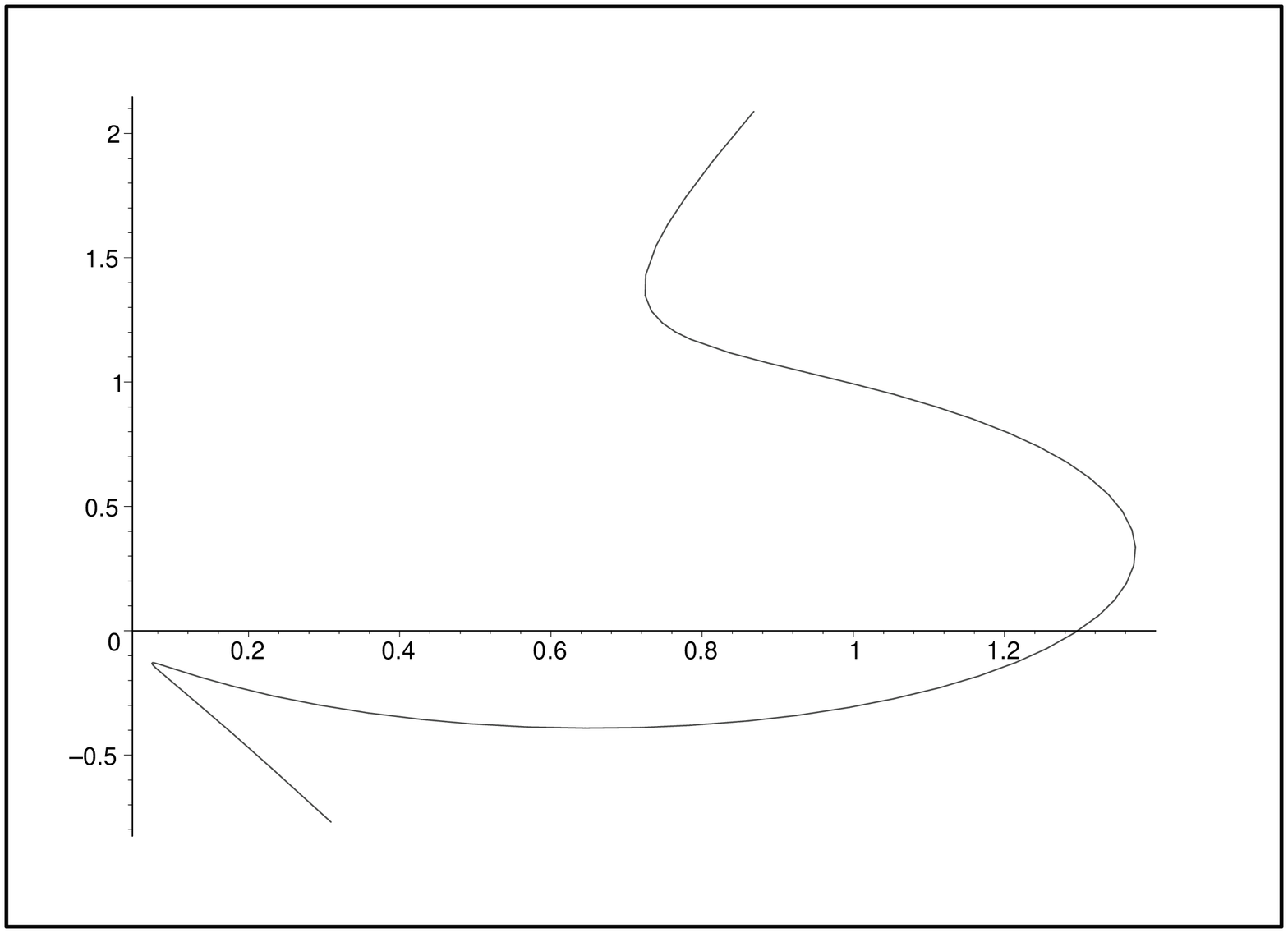}
   \\ \vspace{3mm}
    \includegraphics[angle=270,width=7.7cm,totalheight=9.0cm]{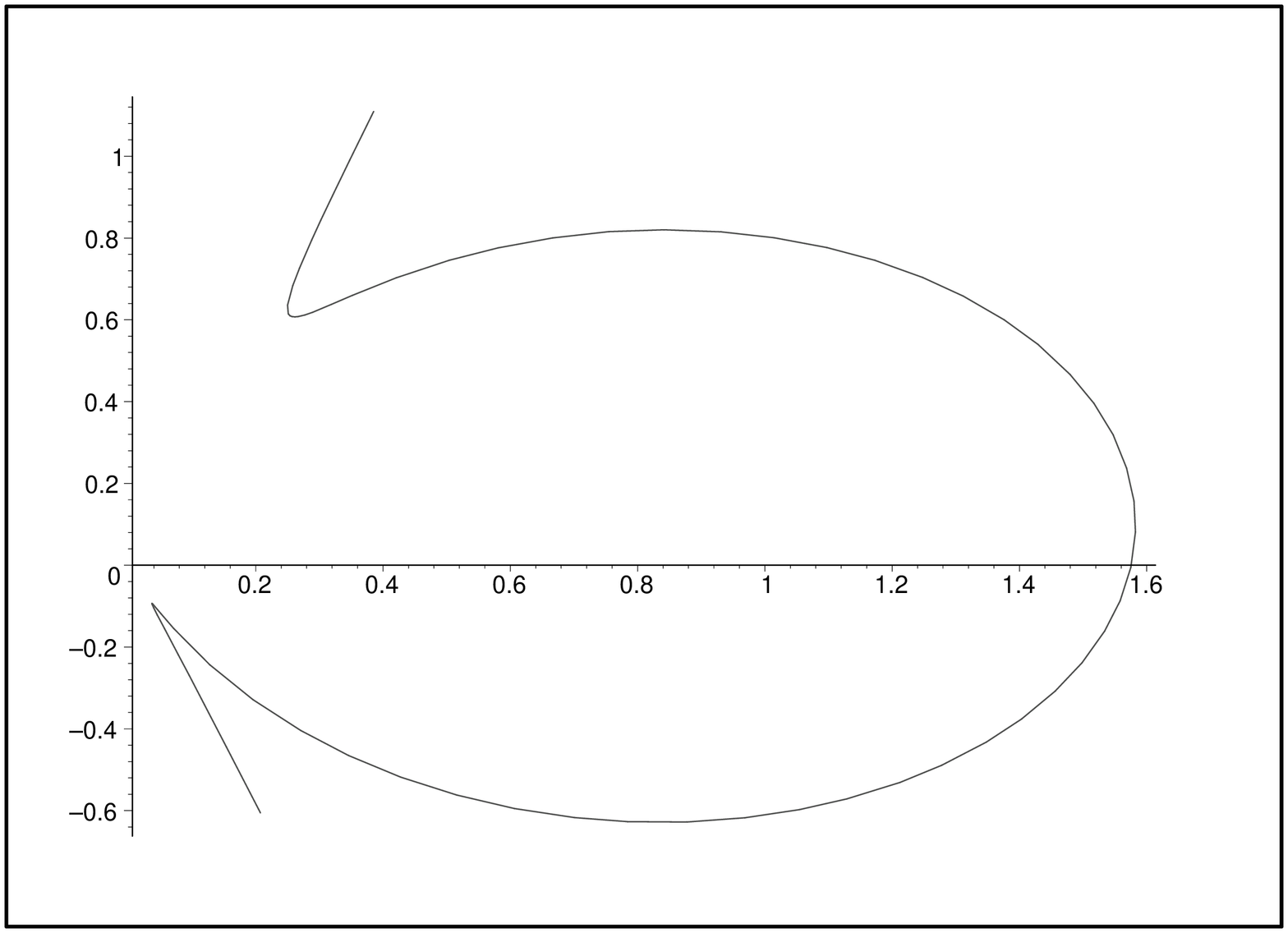}
    \caption{Example 3 (ii), the images of $h_6^{(2)}, h_{10}^{(2)}$ and $h_{30}^{(2)}$.}
\end{figure}

(iii) On the other hand, setting $\mu=1+i$ and choosing
$\tau_n^{(3)}=1-\frac3n$, we find the sequence $h_n^{(3)}$ of
spirallike functions with respect to interior points, which
converges to the function
\[
\widetilde{h}(z)=\left(h(z)\right)^{\frac{1+i}{0.8}}=(1-z)^{1+i}\,,
\]
which is spirallike with respect to a boundary point. In Figure 5,
one can see images of a number of approximating functions and the
image of $\widetilde{h}$.

\vspace{2mm}

\begin{figure}\centering 
    \includegraphics[angle=270,width=6cm,totalheight=6cm]{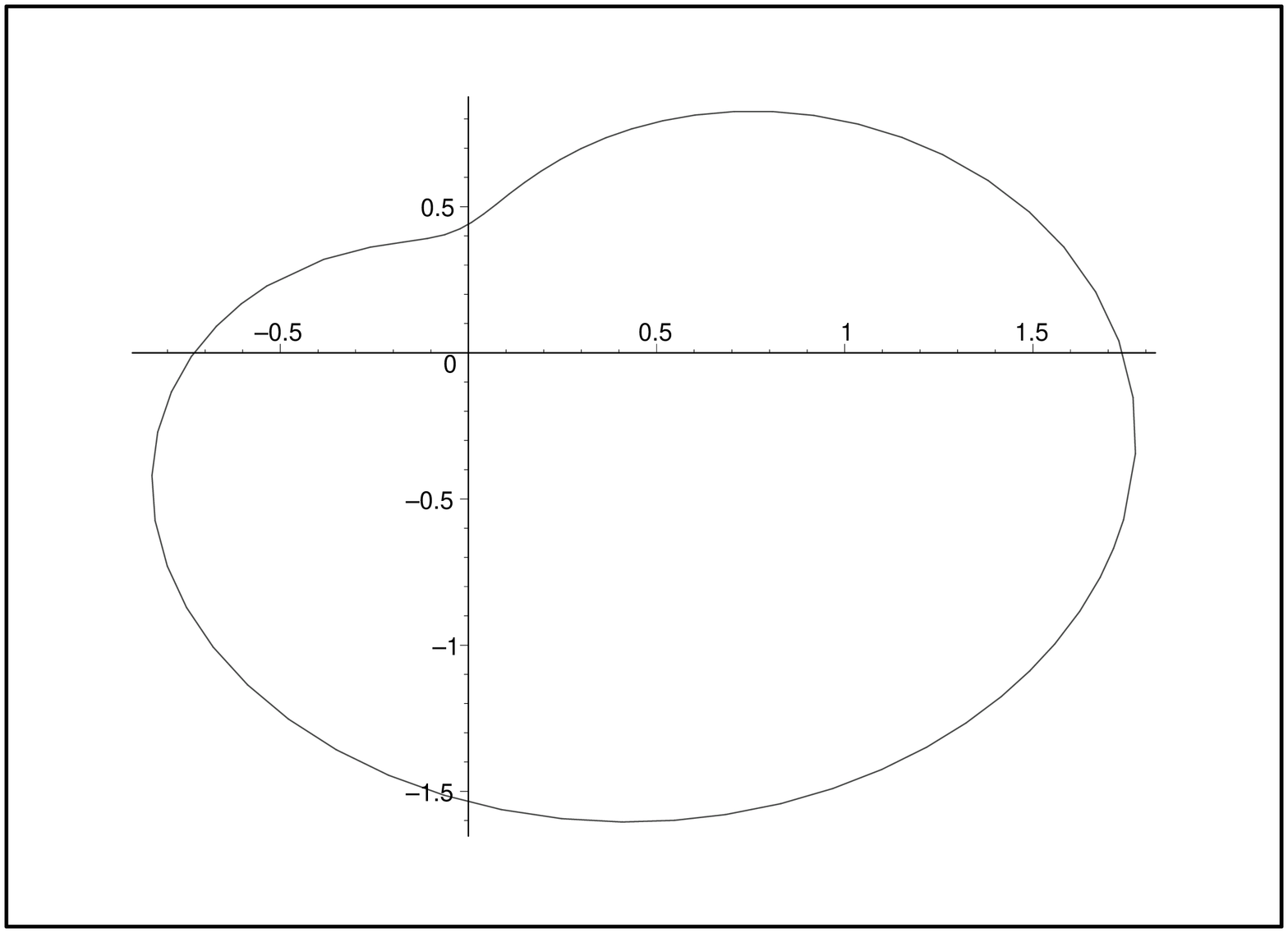}
\hspace{3mm}
    \includegraphics[angle=270,width=6cm,totalheight=6cm]{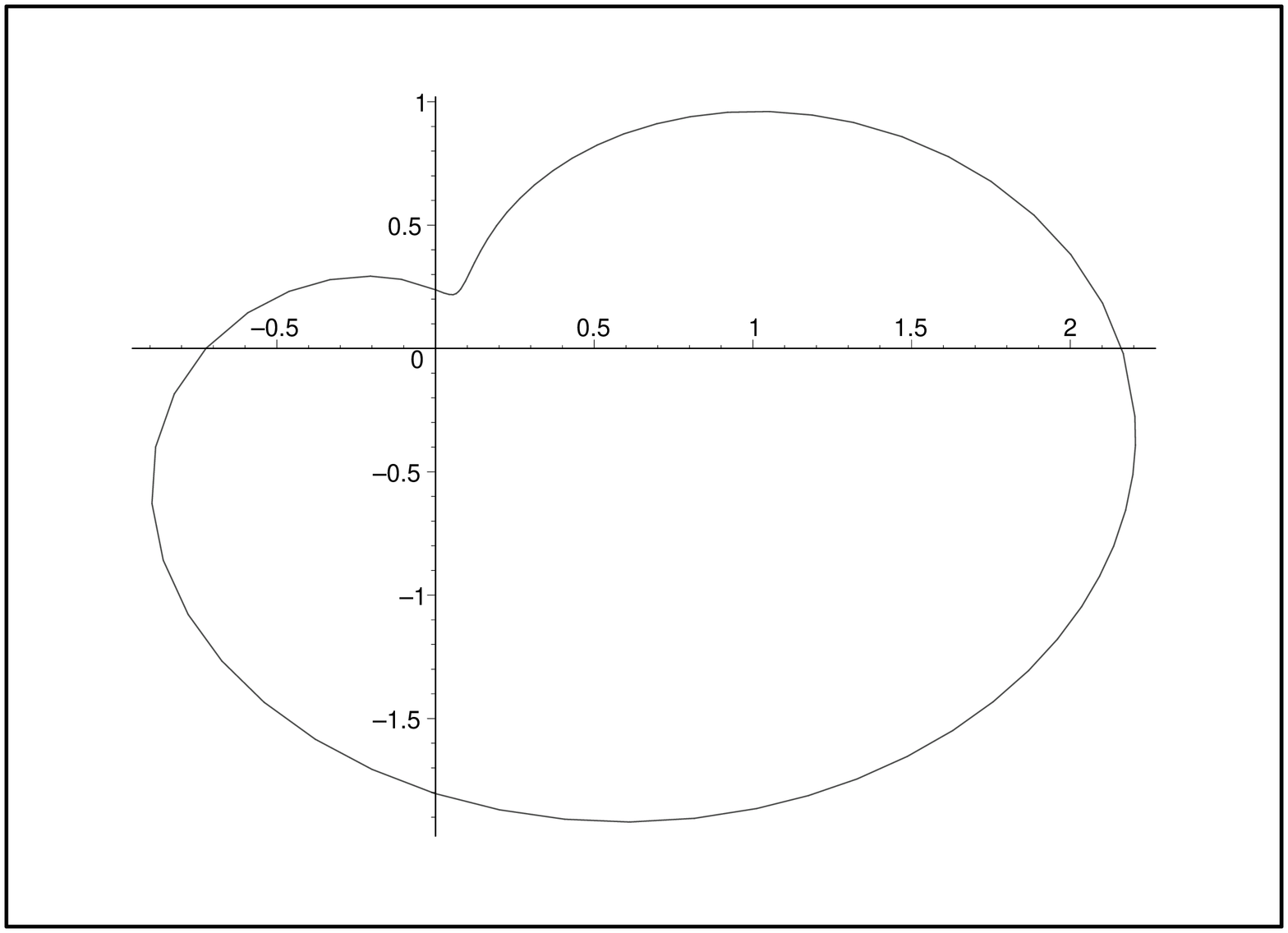}
\\ \vspace{3mm}
    \includegraphics[angle=270,width=6cm,totalheight=6cm]{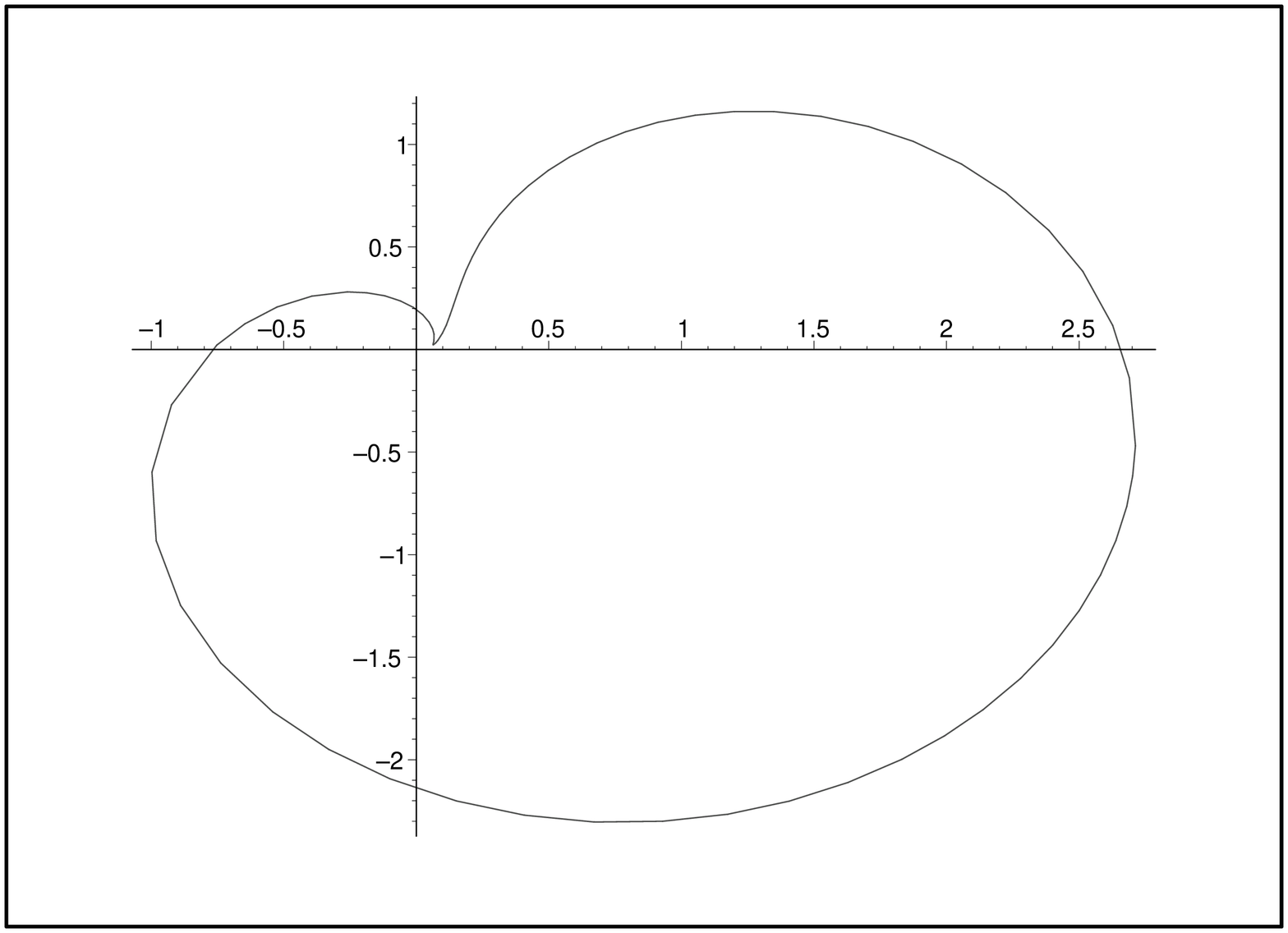}
\hspace{3mm}
    \includegraphics[angle=270,width=6cm,totalheight=6cm]{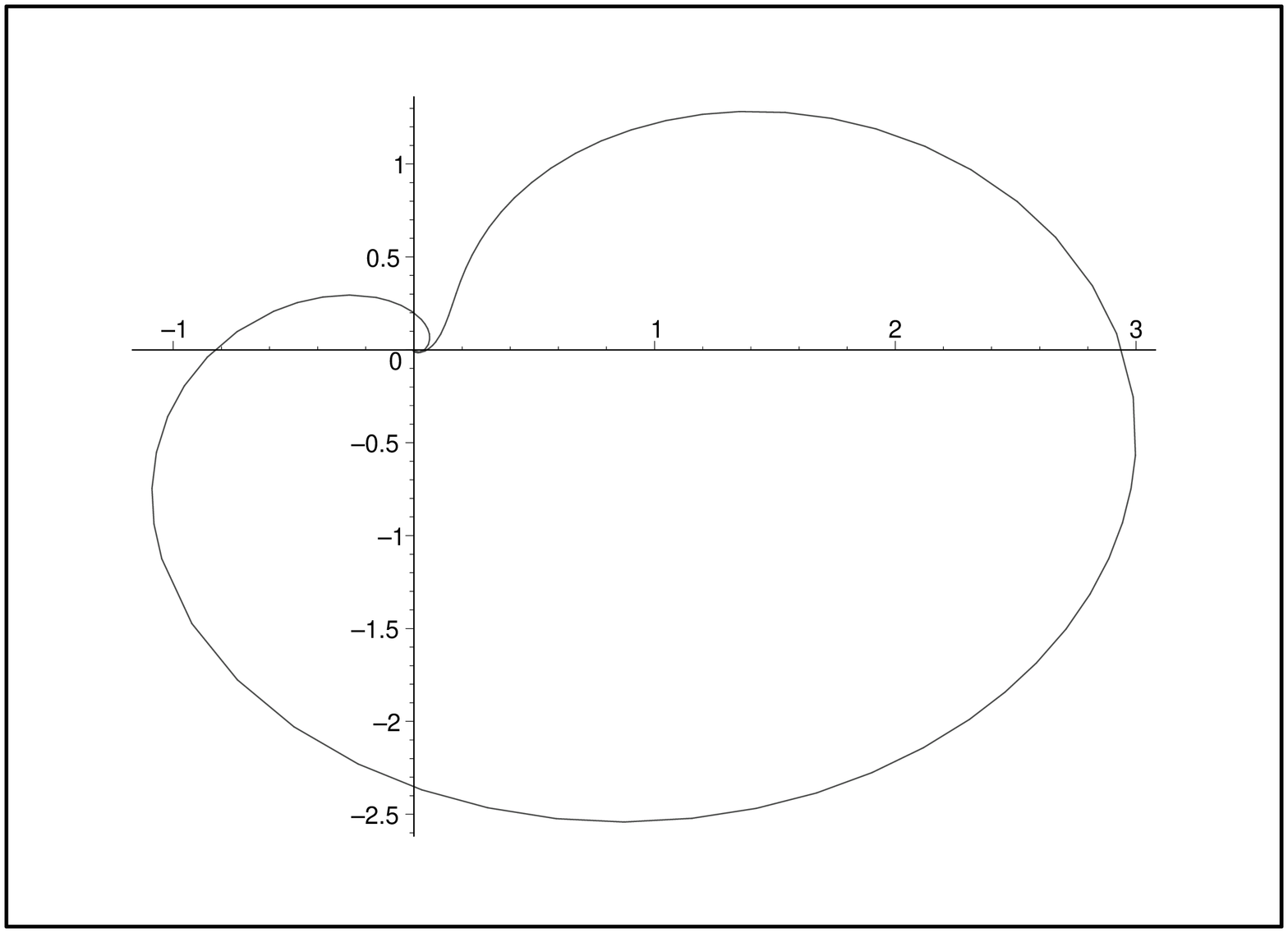}
    \caption{Example 3 (iii), the images of $h_6^{(3)}, h_{10}^{(3)}, h_{30}^{(3)}$ and $\widetilde h$.}
\end{figure}

\vspace{3mm}

Since $\gup$ is a real cone, we can always consider equation
(\ref{26}) with $f$ normalized by $f'(1)=1$.

\begin{defin}
We say that $h\in\Hol(\Delta,\C)$ belongs to the class
$\spi_\lambda[1]$ if it satisfies the equation
\[
\lambda h(z)=h'(z)f(z),
\]
where $f\in\gup$ with $f'(1)=1$ and
$\lambda\in\Omega_+=\left\{w:\, |w-1|\le1,\,w\neq0\right\}$.
\end{defin}

It follows by Theorem 3.4 that
\[
\spi[1]=\bigcup_{\lambda\in\Omega_+}\spi_\lambda[1]
\]
and
\[
S^*[1]=\bigcup_{0<\lambda\le2}\spi_\lambda[1].
\]
In fact, for $\lambda\in(0,2]$, the class $\spi_\lambda[1]$
consists of starlike functions $h$ such that the smallest wedge
which contains the image $h(\Delta)$ is exactly of angle
$\lambda\pi$.

\begin{theorem}
Let $\lambda\in\Omega_+=\left\{w:\ |w-1|\leq1\right\},\
\lambda\ne0.$ Then a function $h\in\Hol(\Delta,\C)$ belongs to the
class $\spi_\lambda[1]$ if and only if it admits the
representation
\be\label{th3.8}
h(z)=(1-z)^{2\lambda}\left[\frac{h_{*}(z)}{z}\right]^\lambda,
\ee
where $h_0$ is a starlike function of class $S^*[0]$ which
satisfies the condition
\[
2\inf\left(\frac{1-|z|^2}{|1-z|^2}\,\Re\frac{z
h_*'(z)}{h_*(z)}\right)=1.
\]
\end{theorem}

\pr Let $h\in\spi[1]$ be $\lambda$-spirallike. Then it satisfies
the equation
\[
\lambda h'(z) =h(z) f(z)\left(=-h(z)
(1-z)^2\cdot\frac1{q(z)}\right),
\]
where $f(z)\in\gup$ with $f'(1)=1$ (or, what is the same,
$q'(1)=-1$ and $q\in\Pp$).

By the above theorem, the functions $h_{\tau}\in\Hol(\Delta,\C)$
defined by
\[
h_\tau(z)
=h^{\frac\mu\lambda}(z)\,\frac{(z-\tau)(1-z\bar\tau)^{\frac\mu{\bar\mu}}}
{-(1-z)^{1+\frac\mu{\bar\mu}}}
\]
are spirallike with respect to an interior point ($h(\tau)=0$) and
satisfy the equations
\[
\mu h_\tau(z)
=h_\tau'(z)(z-\tau)(1-z\bar\tau)\,\frac1{q_\tau(z)}\, ,
\]
where $q_\tau\in\Pp$ and defined by the formula
\[
q_\tau(z)=\frac1z\,\left[(z-\tau)(1-z\bar\tau)r(z)+ \bar\gamma\tau
- \gamma\bar\tau z^2+2iz\Im\gamma\right] ,
\]
where
\[
r(z)=\frac{z q(z)+\gamma z^2-\bar\gamma-2iz\Im\gamma}
{-(1-z)^2}\,,\ \gamma =\frac1\mu\,.
\]

If, in particular, we set $\mu =1$ and $\tau =0$, then we get
$\gamma =1$; and
$h_0(z)=h^{\frac1\lambda}\cdot\frac{z}{-(1-z)^2}$, hence
$h_*=-h_0$, is of the class $S^*[0]$ and satisfies the equation
\[
h_*(z) =h_*'(z) f_0(z) =h_*'(z) \cdot \frac{z}{q_0(z)},
\]
where $q_{0}(z) =r(z).$

Since
\[
\delta_r(1) =\angle\lim_{z\to1}(1-z) r(z)
=\angle\lim_{z\to1}\frac{z
q(z)}{z-1}+\angle\lim_{z\to1}\frac{z^2-1}{z-1}=1,
\]
we obtain
\[
\Re q_0(z) =\Re\frac{zh_*'(z)}{h_*(z)}\geq \frac12\,\frac{1-|z|^2}
{|1-z|^2}.
\]

Conversely, let $h\in\Hol(\Delta,\C)$ admits the representation
\[
h(z) =(1-z)^{2\lambda} \left[\frac{h_*(z)}{z}\right]^\lambda,
\]
where $h_0$ satisfies the required condition. Differentiating
logarithmically, we obtain
\[
\frac1\gamma\,\frac{h'(z)}{h(z)}=
-\frac{2}{1-z}+\frac1z\,\left[\frac{zh_0'(z)}{h_{0}(z)}-1\right]
=-\frac{2}{1-z}+\frac{1}{z}\left[ q_0(z) -1\right] .
\]

Consider the function $q\in\Hol(\Delta,\C)$ defined by
\[
q(z) = -\frac{(1-z)^2}{\lambda}\,\frac{h'(z)}{h(z)}\,.
\]
We have to show that $q\in\Pp$ and $q'(1) =-1.$ This will imply
that $h$ satisfies Definition~3.2 with
$f(z)=-(1-z)^2\frac1{q(z)}$.

Indeed, if $g(z) =z q(z)$, we have
\[
g(z) =1-z^2-(1-z)^2 q_0(z).
\]
Both terms in the right-hand side of this equality are in $\G$.
Hence $g$ also belongs to $\G$, since $\G$ is a real cone.

Thus $q$ must belong to $\Pp$ by the Berkson--Porta formula
(\ref{bp}) with $\tau=0$.

Finally, a direct calculation shows that
\[
q'(1) =\angle \lim_{z\to1}
\frac{q(z)}{z-1}=\angle\lim_{z\to1}\frac{1-z^{2}}{z-1}+\angle\lim_{z\to1}(1-z)q_0(z)=-1,
\]
and we are done.

\vspace{2mm}

Define $\widetilde{S}^*[0]$ as the subclass of $S^*(=S^*[0])$ of
all starlike functions $h\in\Hol(\Delta,\C),\ h(0)=0$, which
satisfy the condition
\[
2\inf \frac{1-|z|^2}{|1-z|^2} \Re\frac{zh'(z)}{h(z)}=1.
\]

\begin{corol}
For each $\lambda\in\Omega_+=\left\{w:\,|w-1|\le1,\ w\neq0
\right\}$, the class $\widetilde{S}^*[0]$ is homeomorphic to the
class $\spi_\lambda[1]$.
\end{corol}

Using representation (\ref{th3.8}), one checks easily the
following characterization of spirallike functions with respect to
a boundary point (see \cite{RMS-81, LA, SH-SEM}).

\begin{corol}[Generalized Robertson condition]
A function \newline ${h\in\Hol(\Delta,\C)}$ normalized by the
condition $h(0)=1$ is of class $\spi[1]$ if and only if for some
$\lambda\in\Omega_+$
\[
\Re\left[ \frac{2z}{\lambda }\frac{h'(z)}{h(z)} +
\frac{1+z}{1-z}\right] \geq 0.
\]
\end{corol}

Using Theorem 3.4 above, one can prove the uniform stability of a
semigroup generated by $f\in\gup$ under any perturbation of its
generator (see Problem in Section 2.3).

Let $\left\{F_{t}\right\}_{t\ge0}$ be a semigroup generated by
$f\in\gup$ and let $f_n\in\G[\tau_n]$, $\tau_n\in\Delta,$ be a
sequence of generators converging to $f$ uniformly on compact
subsets of $\Delta$. If $\left\{F_{t}^{(n)}\right\}_{t\geq 0}$ is
a semigroup generated by $f_n$, then it follows that
$F_t^{(n)}(z)\to\tau_n$ for each $z$ as $t$ tends to $\infty$. In
addition, since $f_n$ and $f$ are locally Lipshitzian on $\Delta,$
the uniqueness property of the Cauchy problem implies that for
each $a\in\Delta$, there exist $r>0$ and $T(=T(r))$ such that for
each $z$ from the circle $|z-a|<r$ and each $t\in[0,T(r)]$, the
sequence of values $F_{t}^{(n)}(z)$ converges to $F_t(z)$ as $n$
tends to infinity. The question is whether this convergence is
uniform on each compact subset of $\Delta\times\R^+$, in other
words, whether for each $0<r<1$ and each $0<T<\infty$ which does
not depend on $r$, the sequence $F_t^{(n)}$ converges to $F_{t}$
uniformly on the set $\overline{\Delta_r}\times[0,T]$, where
$\overline{\Delta_r}=\left\{z\in\C:\ |z|\leq r\right\}$.

We now show that the answer is affirmative.

\begin{theorem}
Let $\left\{F_t\right\}_{t\ge0}$ be a semigroup generated by
$f\in\gup$, and let $f_n\subset\G$ be a sequence of generators
converging to $f$ uniformly on compact subsets of $\Delta$ such
that $f_n(\tau_n)=0$ with $\tau_n\in\Delta$. Assume also that the
closure of the set $\left\{\mu_n=f_n'(\tau_n)\right\}$ does not
contain the origin. Then the sequence
$\left\{F_t^{(n)}(z)\right\}_{t\ge0}$ of semigroups generated by
$f_n$ converges to $F_t(z)$ uniformly on compact subsets of
$\Delta\times\R^+$.
\end{theorem}

\pr By passing to a subsequence, if necessary, we may suppose that
the sequence $\mu_n=f_n'(\tau_n)$ is convergent, say, to
$\mu\in\C$. It follows from our assumption that $\mu\ne0$ and
Theorem~3.3 that $\mu$ must satisfy the condition
$|\mu-\beta|\leq\beta$. Then the sequence $\{h_n\}$ defined by
\[
h_n(z)=\exp \left[\mu\int_0^z \frac{dz}{f_\tau(z)}\right]
\]
converges to the function $h_{1}$ defined by (\ref{31}).

In addition, we have
\begin{equation}
F_{t}^{(n)}(z)=h_n^{-1}\left(e^{-t\mu_n}h_n(z)\right).
\end{equation}
Thus $F_t^{(n)}(z)$ converges uniformly on compact subsets of
$\Delta\times\R^+$ to a function
$G:\,\Delta\times\R^+\mapsto\Delta$ defined by
\begin{equation}\label{34}
G(t,z)=h_1^{-1}\left(e^{-t\mu}h_1(z)\right) .
\end{equation}
On the other hand, we already know that if $h$ is a solution of
the equation $\beta h(z) =h'(z) f(z)$, then
\begin{equation}\label{35}
F_{t}(z) =h^{-1}\left(e^{-t\beta}h(z)\right)
\end{equation}
and
\begin{equation}\label{36}
h_1=h^{\frac\mu\beta}.
\end{equation}
Then we get from (\ref{34}) and (\ref{36})
\[
h^{\frac\mu\beta}\left(G(t,z)\right)=e^{-t\mu}h^{\frac\mu\beta}(z)
\]
or
\[
h\left(G(t,z)\right) =e^{-t\beta}h(z) .
\]
The latter equation and (\ref{35}) imply that $G(t,z)=F_{t}(z)$
for all $z\in\Delta$ and $t\in\R^+$.

\vspace{2mm}

Another useful application of the results given above is the
following description of the spectrum of the infinitesimal
generator of a one-parameter semigroup of composition operators in
the infinite dimensional Frech\'{e}t space $E=\Hol(\Delta,\C)$.

For $f\in\G$, define a linear operator $\Gamma _f$ on $E$ by the
formula
\begin{equation}\label{37}
\left(\Gamma _{f}h\right)(z):=h'(z) f(z)
\end{equation}
(see Section 1.3).

Recall that the spectrum $\sigma\left(\Gamma _{f}\right)$ of the
operator $\Gamma_f$ is the set of all complex numbers
$\lambda\in\C$ for which $\lambda I-\Gamma _{f}$ is not
continuously invertible, where $I$ is the identity operator on
$\Hol(\Delta,\C)$. The point spectrum $\sigma_p\left(
\Gamma_{f}\right)$ of $\Gamma_{f}$ is the subset of
$\sigma\left(\Gamma _{f}\right)$ which consists of its
eigenvalues, i.e.,
\[
\sigma_p\left(\Gamma_f\right)=\left\{\lambda\in\C:\ \left(\lambda
I- \Gamma_f\right) h=0\ \mbox{for some } h\neq 0\right\}.
\]

In this case, a nontrivial solution of the equation $\left(\lambda
I-\Gamma_f\right) h=0$ is called an eigenvector of $\Gamma_f$
corresponding to eigenvalue $\lambda$.

Generally speaking, in infinite dimensional spaces, the spectrum
of a linear operator does not coincide with its point spectrum.
However, since in our situation $f\in\G^{+}[1]$ does not vanish in
$\Delta$, we have $\sigma_{p}(\Gamma _{f})=\sigma(\Gamma_{f})=\C.$

This, for example, implies immediately that the linear semigroup
generated by $\Gamma _{f}$ (and hence, the nonlinear semigroup
$\Ss=\{F_{t}\}_{t\ge0}$ generated by $f$) cannot be
holomorphically extended into a sector containing the positive
real axis (see, for example, \cite[Chapter IX.10]{YK-68}).

Over the past few decades a lot of work address a wide range of
topics concerning properties of composition operators on classical
Banach spaces of analytic functions. In particular, composition
operators have been studied extensively in the setting of Hardy or
Bergman spaces on $\Delta$ (see \cite{CCC-MBD} and \cite{SJH-93}).

In the case of an interior null point of a generator, the
eigenvalues of $\Gamma _{f}$ form a discrete set, while in case of
the boundary Denjoy--Wolff point, the spectrum of $\Gamma_{f}$ may
be larger essentially.

Indeed, as we have already seen in Theorem 3.4, the eigenfunction
$h$ corresponding the eigenvalue $2\beta$, where $\beta=f'(1)$, is
univalent and starlike with respect to a boundary point. Then
$h\in H^p$ for each $p<\frac12$ (see, for example, \cite{Gar}).
Hence, for each positive $\lambda$ the function
$h_1(z):=\left(h(z)\right)^{\frac{2\beta}\lambda}\in H^{q}$ when
$q\lambda<\beta$. Obviously, $h_1$ is an eigenfunction
corresponding to $\lambda$. Therefore, for any Hardy space $H^q$
we have $\sigma(\Gamma_f)\supset\left(0,\frac\beta q\right)$.

This, inter alia, motivates consideration of the univalence of
eigenfunctions of composition operators on the (locally convex)
Fr\'echet space $\Hol(\Delta,\C).$ In a slightly more general
setting, we describe the structure of the spectrum of weighted
composition operators in the context of the $k$-valency of the
corresponding eigenfunctions. Namely, we use the following
definition.

\begin{defin}[see \cite{GAW}, p. 89]
A function $f$ meromorphic in a domain $D$ is said to be
$k$-valent in $D$ if for each $w_0$ (infinity included) the
equation $f(z)=w_0$ has at most $k$ roots in $D$ (where the roots
are counted in accordance with their multiplicity) and if there is
some $w_1$ such that the equation ${f(z)=w_1}$ has exactly $k$
roots in $D$.
\end{defin}

Let $E$ be a space of meromorphic functions in $\Delta,$ and let
$\Ss=\{F_{t}\}_{t\ge0}$ be a semigroup of holomorphic
self-mappings of $\Delta$ generated by $f\in\G^{+}[1],\
f'(1)=\beta.$ For a suitable $w\in E$, one can define a (weighted)
composition semigroup of linear operators $T_{t}:\ E\mapsto E,\
t\geq 0$, by the formula
\[
T_{t}(h)(z)=\frac{w(F_{t}(z))}{w(z)}h(F_{t}(z)).
\]
This semigroup is generated by the operator
$\widetilde{\Gamma}_{f}$ defined by
\be\label{starg}
\widetilde{\Gamma }_{f}h=h^{\prime }f+hf\frac{w^{\prime
}}{w},\quad h\in\Hol(\Delta,\C).
\ee
When $w\equiv 1$, this reduces to the unweighted semigroup of
composition operators $\{C_{t}\}_{t\ge 0}$ generated by
$\Gamma_{f}$.

It is clear that for each $\lambda$ in the point spectrum
$\sigma_{p}(\widetilde{\Gamma }_{f})$, the eigenspace $E_\lambda$
corresponding to $\lambda$ is one-dimensional. For
$k\in\N\bigcup\{\infty\}$, denote by $\sigma^{(k)}$ the subset of
$\sigma_{p}(\widetilde{\Gamma }_{f})$ such that for each
$\lambda\in\sigma^{(k)}$ the function $wh$ is $k$-valent whenever
$h\in E_\lambda$.

\begin{theorem}\label{t3.8}
Let $f\in\gup,\ f'(1)=\beta,$ and let the operator $\Gamma_f:\
E\mapsto E$ be defined by {\rm(\ref{starg})}. Then the spectrum
$\sigma_p\left(\widetilde{\Gamma}_f\right)$ is the whole complex
plane $\C$. Moreover,
\[
\sigma_p\left(\widetilde{\Gamma}_f\right)=
\bigcup\limits_{k\in\NN\cup\{\infty\}}\sigma^{(k)},
\]
where for $k\in\N$
\[
\sigma^{(k)}=\left(k\Omega_+\cup k\Omega_-\right)
\setminus\left((k-1)\Omega_+\cup(k-1)\Omega_-\right)
\]
with $\Omega_\pm=\left\{\omega\not=0:\
|\omega\mp\beta|\le\beta\right\}$, and
\[
\sigma^{(\infty)}=\left\{\lambda\in\C:\ \Re\lambda =0\right\} .
\]

In addition, for each $\lambda\neq 0$ and for each element $h$ of
the eigenspace corresponding to $\lambda$, the following modified
Robertson inequality holds:
\be\label{4ag}
\Re \left[\frac{2\beta}{\lambda}\,\frac{zh'(z)}{h(z)}+
\frac{1+z}{1-z}+\frac{2\beta}{\lambda}\,\frac{zw'(z)}{w(z)}\right]>0.
\ee
\end{theorem}
Fig. 6 illustrates the sets $\sigma^{(k)}$ described in Theorem
3.8.
\begin{figure}\centering 
    \includegraphics[angle=0,width=7cm,totalheight=3.8cm]{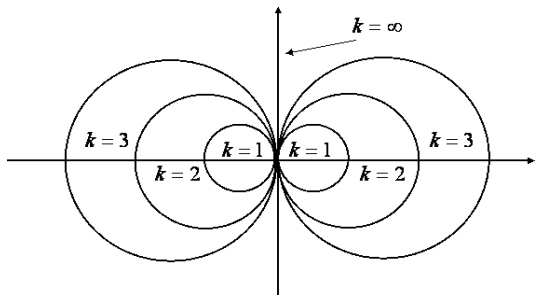}
    \caption{The sets $\sigma^{(k)}$ of $k$-valence of eigenfunctions.}
\end{figure}

\pr The eigenvalue problem for the operator $\widetilde{\Gamma
}_{f}$ is the differential equation
\be\label{1g}
h^{\prime }f+hf\frac{w^{\prime }}{w}=\lambda h.
\ee
To solve it, we denote by $h^*$ the starlike function with respect
to a boundary point which satisfies
\be\label{2g}
\beta h^{\ast }=h^{\ast \prime }f
\ee
and is normalized by $h^*(0)=1,\ h^*(1)=0$ ($h^*$ is the
eigenfunction of $\Gamma _{f}$ corresponding to the eigenvalue
$\beta$).

Then (\ref{1g}) and (\ref{2g}) imply that
\[
\frac{(wh)^{\prime }}{wh}=\frac{\lambda }{\beta
}\,\frac{{h^*}'}{h^*},
\]
and hence
\be\label{3g}
w(z)h(z)=a\left(h^*(z)\right)^{\lambda /\beta },\quad a\in\C.
\ee

For a given eigenvalue $\lambda$, (\ref{3g}) describes the
eigenspace corresponding to $\lambda.$ So it is enough to prove
our assertion for the case $w\equiv1$ or, what is the same,
$\widetilde{\Gamma }_{f}=\Gamma _{f}.$

To this end, let $\lambda\in\C$. First suppose that $\Re\lambda>0$
and \linebreak$\lambda\in k\Omega_+\setminus(k-1)\Omega_+$. Then
${\lambda_1=\frac\lambda k}$ belongs to $\Omega_+$. Thus, by
Theorem~3.4, the solution $h_1$ of the initial value problem
\[
h_1'(z)f(z)=\lambda_1 h_1(z),\quad h_1(0)=1,
\]
is univalent. Set $h(z)=h_1(z)^k$. Evidently, this function
satisfies the equation
\[
\Gamma_f h(z)\left(=h'(z)f(z)\right)=\lambda h(z),
\]
and is at most $k$-valent. To show that $h$ is $k$-valent,
consider the univalent solution $h_2$ of the initial value problem
\[
h_2'(z)f(z)=2\beta h_2(z),\quad h_2(0)=1.
\]
This function is starlike with respect to a boundary point, and
the smallest wedge which contains the image $h_2(\Delta)$ is
exactly of angle $2\pi$ (see the discussion after Definition~3.2).
Therefore, for each $\epsilon>0$, there is real $\psi$ and $r>0$
small enough such that the image $h_2(\Delta)$ contains the curve
\[
\left\{r\exp\left(\phi\,\frac{\Im\lambda}{\Re\lambda}\right)\,\exp\left(i\phi-i\log
r \frac{\Im\lambda}{\Re\lambda} \right):\
\psi\le\phi\le\psi+2\pi(1-\epsilon) \right\}.
\]
Since $\lambda\not\in(k-1)\Omega_+$, one can choose
$\epsilon<1-\frac{2\beta(k-1)\Re\lambda}{|\lambda|^2}\,$. Then
there are $k$ different points $z_0,z_1,\ldots z_{k-1}$ in
$\Delta$ for which
\bep
h_2(z_m)= r\exp\left(\phi_m\,\frac{\Im\lambda}{\Re\lambda}\right)
\,\exp\left(i\phi_m-i\log r \frac{\Im\lambda}{\Re\lambda}
\right),\\ \phi_m=\psi+2\pi
m\,\frac{2\beta\Re\lambda}{|\lambda|^2}\, ,\quad m=0,1,\ldots,k-1.
\eep

Note that $h(z)=h_2(z)^{\lambda/2\beta}$. Now a simple calculation
shows that $h(z_m)$ does not depend on $m=0,1,\ldots,k-1$, i.e.,
$h$ is $k$-valent and $\lambda\in\sigma^{(k)}$.

The case $\Re \lambda<0$ is reduced to the previous one by
replacing $h$ by $1/h$.

Finally, suppose $\Re\lambda=0$. As usual, denote by
$\left\{F_t\right\}_{t\ge0}$ the semigroup generated by $f$. Thus,
the function $u(t,z):=F_t(z)$ satisfies the Cauchy problem
\bep
\left\{
\begin{array}{l}
\dst\frac{\partial u(t,z)}{\partial t}+f(u(t,z))=0\\[3mm]
u(0,z)=z\in\Delta.
\end{array}
\right.
\eep
Solving this problem for $z=0$, we have
\[
\int_0^{u(t,0)}\frac{du}{f(u)}=-\int_0^tdt=-t.
\]
Since the function $1/f$ is holomorphic on the open disk $\Delta$,
we conclude that for each point $z$ which belongs to the curve
$\Lambda:=\left\{u(t,0):\, t\ge0 \right\}$ (joining the origin
with the boundary point $z=1)$, the integral
\[
\int_0^z\frac{du}{f(u)}
\]
takes real values which tend to $-\infty$ as $z\in\Lambda$ tends
to $1$.

Returning to the eigenvalue problem
\[
\Gamma_f(h)(z)\left(:=h'(z)f(z)\right)=\lambda h(z)
\]
and separating variables in this differential equation, we have
\[
\int_0^z\frac{dh(z)}{h(z)}=\int_0^z\frac{\lambda dz}{f(z)}\,,
\]
so that
\[
\log\left(\frac{h(z)}{h(0)}\right)=\lambda\int_0^z\frac{dz}{f(z)}\,.
\]
Thus, for each point $z\in\Lambda$,
\[
h(z)=h(0)\exp\left(-\lambda t\right),
\]
where $z=u(t,0)$. We claim that this function is infinite-valent.
When $\lambda=0$, this is evident. Suppose $\lambda\not=0$; then
$\lambda=ia,\ a\not=0$. Choosing the sequence $\dst
z_n=u\left(t_0+\frac{2\pi n}a,0\right)\in\Lambda\subset\Delta$ for
any fixed $t_0\ge0$, we see immediately that the value $h(z_n)$
does not depend on $n$.

Now, let $h^*$ be the solution of the differential equation
(\ref{2g}) normalized by $h^*(0)=1,\ h^*(1)=0$. As in the proof of
Theorem~3.3, we conclude that the image of the function
$h^*(\Delta)$ must lie in a wedge of angle $\pi$. Since $h^*$ is
starlike with respect to a boundary point, it follows from a
result of Lyzzaik \cite{LA} (see also \cite{RMS-81}) that
\[
\Re\left[ \frac{2z{h^*}'(z)}{h^*} +\frac{1+z}{1-z} \right]>0.
\]
Substituting (\ref{3g}) into the latter inequality, we get
(\ref{4ag}). The proof is complete.

\vspace{2mm}

Combining this assertion with previous results, one can obtain
additional geometrical information on (characteristics of)
eigenfunctions.

\vspace{3mm}

\noindent\textbf{Example 4.} For each $a\in(0,1)$, we rewrite
(\ref{4ag}) in the form
\be\label{5g}
\Re\left[\frac{2a\beta}\lambda\,\frac{zh'(z)}{h(z)}+\frac{1+z}{1-z}\right]>
\Re\left[(1-a)\,\frac{1+z}{1-z}-
\frac{2a\beta}\lambda\,\frac{zw'(z)}{w(z)} \right].
\ee

Suppose that $w(z)=(1-z)^\delta$ for some $\delta\in\C$.

Let $c$ be a positive number less than $1$. In this case, for each
$\lambda$ such that $\frac\delta\lambda$ is real with
$\frac\delta\lambda<-\frac{c}{a\beta}$, a simple calculation shows
that the right-hand side in (\ref{5g}) is greater than $c$.
Indeed, for each $z\in\partial\Delta$,
\bep
\Re\left[(1-a)\,\frac{1+z}{1-z}-
\frac{2a\beta}\lambda\,\frac{zw'(z)}{w(z)} \right] =
\frac1{|1-z|^2}\,\Re\left[-\frac{2a\beta\delta}\lambda
+\frac{2a\beta\delta z}\lambda \right]\\
=-\frac{2a\beta\delta}\lambda\,\Re\frac{1-z}{|1-z|^2}\ge-\frac{a\beta\delta}\lambda>c.
\eep
Letting $a\to1^-$, we conclude that
\[
\Re\left[\frac{2\beta}\lambda\,\frac{zh'(z)}{h(z)}+\frac{1+z}{1-z}\right]>c.
\]

Therefore, by Theorem 3.1 in \cite{E-pr}, the image $h(\Delta)$ of
the eigenfunction $h$ covers the image of the function
$(1-z)^{\frac{c\lambda}\beta}$.

\vspace{3mm}

\noindent\textbf{Example 5.} Consider now the (nonanalytic) weight
$w(z)=\frac{(1-z)^2}z\,$. By formula~(\ref{3g}), each
eigenfunction has the form
\[
h(z)=a\,\frac{z}{(1-z)^2}\left(h^*(z)\right)^{\frac\lambda\beta}.
\]
In this situation, for all positive $\lambda$ less than $2\beta$,
inequality (\ref{4ag}) becomes
\[
\Re\frac{zh'(z)}{h(z)}>\left(1-\frac\lambda{2\beta}\right)\,\frac{1-|z|^2}{|1-z|^2}>0,
\]
i.e., all eigenfunctions are starlike with respect to an interior
point ($h(0)=0$).

In particular, the function
\[
h_0(z)=\frac{z}{(1-z)^2}\,h^*(z)
\]
is the eigenfunction corresponding the eigenvalue $\lambda=\beta$.
Since $h^*\in\spi_1[1]$, Theorem~3.6 implies that $h_0$ is a
starlike function of class $S^*[0]$ which satisfies the condition
\[
2\inf\left(\frac{1-|z|^2}{|1-z|^2}\,\Re\frac{z
h_0'(z)}{h_0(z)}\right)=1.
\]

\vspace{3mm}

Finally, we observe that the above theorem implies the following
nice characterization of multivalent starlike functions with
respect to a boundary point suggested by D.~Bshouty and A.~Lyzzaik
(see \cite{B-L1} and \cite{B-L2}; cf. also \cite{LeA} and
\cite{LA-LA}).

\begin{theorem}
Let $h\in\Hol(\Delta,\C)$ satisfy the equation
\be\label{ag}
-(1-z)^2\,\frac{h'(z)}{h(z)}=4\,\frac{1-\omega(z)}{1+\omega(z)}\,,
\ee
where $\omega\in\Hol(\Delta)$ is a holomorphic self-mapping of
$\Delta$ with the boundary fixed point $\tau=1$ and positive
multiplier
\[
\alpha=\angle\lim_{z\to1}\frac{1-\omega(z)}{1-z}\,.
\]
Then

(i) $h(\Delta)$ is a starlike domain;

(ii) $h$ is a $k$-valent function if and only if $k-1<\alpha\le
k$.

In particular, if $\tau=1$ is the Denjoy--Wolff point for $\omega$
(whence\linebreak $0<\alpha\le1$), then $h$ is univalent and the
smallest wedge which contains $h(\Delta)$ is of angle
$2\alpha\pi$.
\end{theorem}

\pr Let
\[
p(z)=\frac\alpha2\,\frac{1+\omega(z)}{1-\omega(z)}\,.
\]
This function has positive real part. Moreover,
\[
\angle\lim_{z\to1}(1-z)p(z)=\angle\lim_{z\to1}\frac\alpha2\,
\frac{1-z}{1-\omega(z)}\,(1+\omega(z))=1
\]
Define a generator $f\in\G(\Delta)$ by the Berkson--Porta formula
\[
f(z)=-(1-z)^2p(z).
\]

It is clear that $f\in\gup$ with $\beta=f'(1)=1$.

Since $h$ satisfies the equation
\[
2\alpha h(z)=h'(z)f(z),
\]
and $\alpha$ is positive, applying Theorem~3.8 with
$\lambda=2\alpha$ and $\beta=1$, we obtain the result.

If $\tau=1$ is the Denjoy--Wolff point for $\omega$, i.e.,
$0<\alpha\le1$, then
\[
\angle\lim_{z\to1}\frac{(z-1)h'(z)}{h(z)}=2\alpha\le2;
\]
and the assertion follows by Theorem~7 in \cite{E-R-S2001a}. The
proof is complete.

\end{document}